\documentclass{article}

\usepackage[a4paper]{geometry}

\usepackage[utf8]{inputenc}
\usepackage[T1]{fontenc}
\usepackage{lmodern}

\usepackage{amsmath}
\usepackage{amssymb}
\usepackage{theorem}
\usepackage{enumerate}

\usepackage[pdftex,
            pdfauthor={Jethro van Ekeren, Sven Möller and Nils R. Scheithauer},
            pdftitle={Construction and classification of holomorphic vertex operator algebras},
            pdfsubject={},
            pdfkeywords={},
            pdfproducer={},
            pdfcreator={}]{hyperref}

\newcommand{\al}{\alpha}
\newcommand{\bt}{\beta}

\newcommand{\C}{\mathbb C}

\newcommand{\Q}{\mathbb Q}
\newcommand{\R}{\mathbb R}
\newcommand{\Z}{\mathbb Z}
\newcommand{\la}{\langle}
\newcommand{\ra}{\rangle}

\newcommand{\Aut}{\operatorname{Aut}}
\newcommand{\End}{\operatorname{End}}

\newcommand{\sign}{\operatorname{sign}}

\newcommand{\orb}{\operatorname{orb}}
\newcommand{\tr}{\operatorname{tr}}
\newcommand{\vac}{{\bf 1}}
\newcommand{\ch}{\operatorname{ch}}
\newcommand{\id}{\operatorname{id}}

\newcommand{\wt}{\operatorname{wt}}
\newcommand{\eop}{\hspace*{\fill} $\Box$}

\renewcommand{\S}{\mathcal{S}}

\newcommand{\Irr}{\operatorname{Irr}}
\newcommand{\SL}{\operatorname{SL}}
\newcommand{\GL}{\operatorname{GL}}
\newcommand{\Or}{\operatorname{O}}
\newcommand{\V}{\mathcal{V}}

\theoremstyle{break}

\newtheorem{thm}{Theorem}[section]
\newtheorem{prp}[thm]{Proposition}
\newtheorem{cor}[thm]{Corollary}

\begin{document}

\begin{center}
{\Large \bf Construction and classification of \\[2mm]
holomorphic vertex operator algebras} \\[12mm]
Jethro van Ekeren\textsuperscript{1},
Sven Möller\textsuperscript{2},
Nils R.\ Scheithauer\textsuperscript{2}\\[4mm]
\textsuperscript{1}{Universidade Federal Fluminense, Niterói, RJ, Brazil}\\
\textsuperscript{2}{Technische Universität Darmstadt, Darmstadt, Germany}
\end{center}

\vspace*{10mm}
\noindent
We develop an orbifold theory for finite, cyclic groups acting on holomorphic vertex operator algebras. Then we show that Schellekens' classification of $V_1$-structures of meromorphic conformal field theories of central charge $24$ is a theorem on vertex operator algebras. Finally we use these results to construct some new holomorphic vertex operator algebras of central charge $24$ as lattice orbifolds.
\vspace*{10mm}

\section{Introduction}

Vertex algebras give a mathematically rigorous description of 2-dimensional quantum field theories. They were introduced into mathematics by R.\ E.\ Borcherds \cite{B1}. The most prominent example is the moonshine module $V^{\natural}$ constructed by Frenkel, Lepowsky and Meurman \cite{FLM}. This vertex algebra is $\Z$-graded and carries an action of the largest sporadic simple group, the monster. Borcherds showed that the corresponding twisted traces are hauptmoduls for genus $0$ groups \cite{B2}. 

The moonshine module was the first example of an orbifold in the theory of vertex algebras. 
The $(-1)$-involution of the Leech lattice $\Lambda$ lifts to an involution $g$ of the associated lattice vertex algebra $V_{\Lambda}$. The fixed points $V_{\Lambda}^g$ of $V_{\Lambda}$ under $g$ form a simple subalgebra of $V_{\Lambda}$. Let $V_{\Lambda}(g)$ be the unique irreducible $g$-twisted module of $V_{\Lambda}$ and $V_{\Lambda}(g)_{\Z}$ the subspace of vectors with integral $L_0$-weight. Then the vertex algebra structure on $V_{\Lambda}^g$ can be extended to $V_{\Lambda}^g \oplus V_{\Lambda}(g)_{\Z}$. This sum is the moonshine module $V^{\natural}$. In contrast to $V_{\Lambda}^g$, which has 4 irreducible modules, the moonshine module $V^{\natural}$ has only itself as irreducible module.

It was believed for a long time that the above construction can be generalised to automorphisms of any finite order (see. e.g. \cite{Don94}). We show that under certain regularity assumptions this conjecture is true. More precisely we prove (see.\ Theorem \ref{orbtheorem}):

{\em Let $V$ be a simple, rational, $C_2$-cofinite, holomorphic vertex operator algebra of CFT-type and let $g$ be an automorphism of $V$ of order $n$ and type $0$. Then the cyclic group $G = \langle g \rangle$ acts naturally on the irreducible $g^i$-twisted $V$-modules $V(g^i)$ and these actions are unique up to multiplication by an $n$-th root of unity. Suppose the modules $V(g^i)$ have positive conformal weight, i.e. positive smallest $L_0$-eigenvalue, for $i \neq 0 \! \mod n$. Then the actions can be chosen such that   
\[  V^{\orb(G)} = \bigoplus_{i \in \Z_n} V(g^i)^G   \]
has the structure of a vertex operator algebra extending that of $V^G$. With this structure $V^{\orb(G)}$ is simple, rational, $C_2$-cofinite, holomorphic and of CFT-type.} 

We sketch the proof. By \cite{DM0}, \cite{M2} and \cite{CM} the vertex operator algebra $V^G$ is simple, rational, $C_2$-cofinite, self-contragredient and of CFT-type. The irreducible modules of $V^G$ are given up to isomorphism by the eigenspaces $W^{(i,j)} = \{ v \in V(g^i) \, | \, g.v = e(j/n) v \}$
\cite{MT}. Using the Ver\-linde formula proved by Huang \cite{H3} for vertex operator algebras we show that the modules $W^{(i,j)}$ are simple currents and that the fusion group is isomorphic to $\Z_n^2$.
This implies that the sum $W = \bigoplus_{(i,j) \in \Z_n^2} W^{(i,j)}$ is an abelian intertwining algebra whose braiding is determined by the conformal weights of the modules $W^{(i,j)}$. From this the theorem can be derived.

We remark that our approach is inspired by Miyamoto's theory of $\Z_3$-orbifolds of lattice vertex algebras \cite{M1}.

Another important problem in the theory of vertex algebras is the classification of vertex operator algebras satisfying the above regularity conditions.

The rank of a positive definite, even, unimodular lattice is divisible by $8$. In rank $8$ there is exactly one such lattice up to isomorphism, the $E_8$-lattice. In rank $16$ there are two isomorphism classes, $E_8^2$ and $D_{16}^+$. Niemeier showed that in rank $24$ there are exactly $24$ isomorphism classes, the Niemeier lattices, and that the isomorphism class of such a lattice is determined by its root sublattice, the sublattice generated by the elements of norm $2$. The Minkowski-Siegel mass formula shows that with increasing rank the number of classes grows very rapidly. For example in rank $32$ there are more than $10^7$ classes.  

A similar result is expected to hold for simple, rational, $C_2$-cofinite, holomorphic vertex operator algebras of CFT-type. Examples of such vertex operator algebras are obtained from positive definite, even, unimodular lattices. The central charge of such a vertex operator algebra is a positive integer multiple of $8$. Dong and Mason \cite{DM1} proved that for central charge $8$ there is exactly one isomorphism class, the vertex operator algebra $V_{E_8}$ of the $E_8$-lattice, and for central charge $16$ there are two classes, $V_{E_8^2}$ and $V_{D_{16}^+}$. The vectors of weight $1$ in a vertex operator algebra of CFT-type form a Lie algebra. We show that the classification of $V_1$-structures given by Schellekens is a theorem on vertex operator algebras (Theorem \ref{Thm.schellekens}):

{\em Let $V$ be a simple, rational, $C_2$-co\-finite, holomorphic vertex operator algebra of CFT-type and central charge $24$. Then either $V_1 = 0$ or $ \dim (V_1) = 24$ and $V$ is isomorphic to the vertex operator algebra of the Leech lattice or $V_1$ is one of $69$ semisimple Lie algebras described in Table $1$ of \cite{S}.} 

The argument is essentially the one given by Schellekens. The vertex operator subalgebra $W$ of $V$ generated by the elements in $V_1$ is isomorphic to a tensor product of affine vertex operator algebras. It is rational with well-known irreducible modules \cite{FZ}. The vertex operator algebra $V$ decomposes into finitely many irreducible $W$-modules under the action of $W$. The regularity conditions imply that the character of $V$ is a Jacobi form of weight $0$ and index $1$ \cite{KM}. This imposes strong restrictions on $W$ and the structure of $V$ as $W$-module. In particular the possible isomorphism types of $W$ are completely fixed and given in Schellekens list. The main difference to \cite{S} is that we write down explicitly the equations which imply the classification in terms of homogeneous weight polynomials.

It is believed that all the Lie algebras given in Schellekens list are realised as $V_1$-structures of holomorphic vertex operator algebras of central charge $24$ and that the $V_1$-structure determines the vertex operator algebra up to isomorphism. 

Using the main results described above we construct $5$ new holomorphic vertex operator algebras of central charge $24$ as orbifolds of Niemeier lattices and determine their $V_1$-structure (Theorem \ref{NewVOATheorem}):
 
{\em There exist holomorphic vertex operator algebras of central charge $24$ with the Lie algebra $V_1$ given by $A_{2,1} B_{2,1} E_{6, 4}$, $A_{4,5}^2$, $A_{2,6} D_{4,12}$, $A_{1,1} C_{5,3} G_{2,2}$ and $C_{4, 10}$.}

Together with recent results by Lam et al.\ \cite{LS3}, \cite{LS4}, \cite{LL} this implies that the first part of the above classification conjecture holds.

Another application of our results is Carnahan's proof of Norton's generalised moonshine conjecture \cite{C}.

The paper is organised as follows.

First we recall some important modularity properties of vertex operator algebras.

Let $V$ be a simple, rational, $C_2$-cofinite, self-contragredient vertex operator algebra of CFT-type. Suppose that all irreducible $V$-modules are simple currents so that the fusion algebra of $V$ is the group algebra of some finite abelian group $D$ and that the irreducible modules different from $V$ all have positive con\-formal weight.

In section $3$ we show that the conformal weights of the irreducible modules of $V$ define a quadratic form $q_{\Delta}$ on $D$ and that Zhu's representation is up to a character the Weil representation of $D$ (Theorem \ref{thm.3.4}). 

The direct sum of the irreducible $V$-modules has the structure of an abelian intertwining algebra whose associated quadratic form is $-q_{\Delta}$ (Theorem \ref{thm.4.7}). Restricting the sum to an isotropic subgroup gives a vertex operator algebra extending $V$ (Theorem \ref{extensionthm}).

Next assume that $V$ is a simple, rational, $C_2$-cofinite, holomorphic vertex operator algebra of CFT-type and let $G = \langle g \rangle$ be a finite, cyclic group of order $n$ acting by auto\-mor\-phisms on $V$. Suppose that the $g^i$-twisted modules $V(g^i)$ of $V$ have positive conformal weight for $i \neq 0$.

In section $5$ we show that the irreducible modules of $V^G$ are simple currents. This implies that the fusion algebra of $V^G$ is the group algebra of some finite abelian group $D$. We show that $D$ is a central extension of $\Z_n$ by $\Z_n$ whose iso\-mor\-phism type is determined by the conformal weight of the twisted module $V(g)$ (Theorem \ref{thm.5.12}). Furthermore the direct sum of all irreducible $V^G$-modules is an abelian intertwining algebra with the braiding determined by the conformal weights of the irreducible $V^G$-modules and the restriction of the sum to an isotropic subgroup $H$ of $D$ is a vertex operator algebra extending $V^G$ called the orbifold of $V$ with respect to $G$ and $H$ (Theorem \ref{orbtheorem}).

In section $6$ we show that Schellekens' classification of $V_1$-structures of meromorphic conformal field theories of central charge $24$ is a rigorous theorem on vertex operator algebras. 

Next we recall some results on lattice vertex algebras.

In the last section we apply our results to construct $5$ new holomorphic vertex operator algebras of central charge $24$ as orbifolds of Niemeier lattices. We show that they have $V_1$-structures $A_{2,1} B_{2,1} E_{6, 4}$, $A_{4,5}^2$, $A_{2,6} D_{4,12}$, $A_{1,1} C_{5,3} G_{2,2}$ and $C_{4, 10}$ (Theorem \ref{NewVOATheorem}).

We thank P.\ Bantay, S.\ Carnahan, T.\ Creutzig, G.\ Höhn, \mbox{Y.-Z.} Huang, V.\ Kac, C.-H.\ Lam, M.\ Miyamoto, A.\ Schellekens, H. Shimakura and H.\ Yamauchi for valuable discussions and the referee for suggesting several improvements.

The first author was supported by a grant from the Alexander von Humboldt Foundation and later by CAPES-Brazil. The second author was partially supported by a scholarship of the Studienstiftung des deutschen Volkes. The second and third author both were supported by the DFG-project ``Infinite-dimensional Lie algebras in string theory''.

After completion of this work we learned of \cite{DRX} which has some overlap with our Section \ref{Section.orb}, in particular Proposition \ref{prop.mods.are.simple}.

\section{Modular invariance and the Verlinde formula}

In this section we recall some important modularity properties of vertex operator algebras.

\vspace*{2mm}
Let $V$ be a rational, $C_2$-cofinite vertex operator algebra of central charge $c$ and $W$ an irreducible $V$-module. We denote the conformal weight of $W$, 
i.e. the smallest $L_0$-eigenvalue of $W$, by $\rho$. Then $c$ and $\rho$ are rational numbers (\cite{DLM}, Theorem 11.3). For $v \in V$ we define the formal sum 
\[  T_W(v,q) = \tr_W o(v)q^{L_0-c/24} = q^{\rho-c/24} \sum_{n=0}^\infty \tr_{W_{\rho+n}} o(v) q^n
\]
where $o(v) = v_{\wt(v)-1}$ for homogeneous $v$, extended linearly to $V$. Zhu has shown (\cite{Z}, Theorems 5.3.2 and 4.2.1)

\begin{thm}
Let $V$ be a rational, $C_2$-cofinite vertex operator algebra of central charge $c$ and $W \in \Irr(V)$, the set of isomorphism classes of irreducible $V$-modules. Then
\begin{enumerate}[i)]
\item 
Let $q = e^{2 \pi i \tau}$. Then the formal sum $T_W(v,\tau)$ converges to a holomorphic function on the complex upper halfplane $H$. 
\item
Let $v \in V_{[k]}$ be of weight $k$ with respect to Zhu's second grading. Then there is a representation 
\[  \rho_V :  \SL_2(\Z) \to \GL(\mathcal{V}(V))  \]
of $\SL_2(\Z)$ on the fusion algebra $\mathcal{V}(V) = \bigoplus_{W \in \Irr(V)} \C W$ of $V$ such that
\[  T_W(v,\gamma \tau) = (c\tau+d)^k \sum_{M \in \Irr(V)} \rho_V(\gamma)_{W,M} T_M(v, \tau)  \]
for all $\gamma = \left(\begin{smallmatrix} a & b \\ c & d \end{smallmatrix} \right) \in \SL_2(\Z)$.
\end{enumerate}
\end{thm}

We denote the images of the standard generators $S = \left( \begin{smallmatrix} 0 & -1 \\ 1 & 0 \end{smallmatrix} \right)$ and $T = \left( \begin{smallmatrix} 1 &  1 \\ 0 & 1 \end{smallmatrix} \right)$ of $\SL_2(\Z)$ under Zhu's representation $\rho_V$ by $\mathcal{S}$ and $\mathcal{T}$. The definition of $T_M(v,\tau)$ as trace implies
\[  \mathcal{T}_{M,N} = \delta_{M,N} e((\rho(M)-c/24))   \]
where $e(x)$ denotes $e^{2\pi i x}$. The $\mathcal{S}$-matrix is related to the fusion coefficients by the Verlinde formula proved for vertex operator algebras by Huang (\cite{H3}, section~5).

\begin{thm}
Let $V$ be a simple, rational, $C_2$-cofinite, self-contragredient vertex operator algebra of CFT-type. Then 
\begin{enumerate}[i)]
\item
The matrix $\mathcal{S}$ is symmetric and $\mathcal{S}^2$ is the permutation matrix sending $M$ to its contragredient module $M'$. Moreover $\mathcal{S}_{V,U} \neq 0$ for $U \in \Irr(V)$.
\item 
The fusion coefficients are given by
\[ \mathcal{N}_{M,N}^{W} = 
\sum_{U \in \Irr(V)} \frac{\mathcal{S}_{M,U} \mathcal{S}_{N,U} \mathcal{S}_{W',U}}{\mathcal{S}_{V,U}}  \, .  \]
\end{enumerate}
\end{thm} 

Now let $V$ be a simple, rational, $C_2$-cofinite vertex operator algebra of central charge $c$ which is holomorphic, i.e. $V$ has only one irreducible module, and let $G = \langle g \rangle$ be a finite, cyclic group of automorphisms of $V$ of order $n$. Then for each $h \in G$ there is a unique irreducible $h$-twisted $V$-module $V(h)$ of conformal weight $\rho(V(h)) \in \Q$ and a representation 
\[  \phi_h:G \rightarrow \Aut_{\C}(V(h))  \]
of $G$ on the vector space $V(h)$ such that
\[  \phi_h(k)Y_{V(h)}(v,z)\phi_h(k)^{-1} = Y_{V(h)}(kv,z)  \]
for all $k \in G$ and $v\in V$ \cite{DLM}. By Schur's lemma this representation is unique up to multiplication by an $n$-th root of unity. Setting $h = g^i$, $k = g^j$ we define the twisted trace functions
\[  T(v,i,j, q) 
= \tr_{V(g^i)} o(v) \phi_{g^i}(g^j)q^{L_0-c/24}   \, . \]
Dong, Li and Mason have shown (\cite{DLM}, Theorem 1.4)

\begin{thm}\label{DLM.thm}
Let $V$ be a simple, rational, $C_2$-cofinite, holomorphic vertex operator algebra of central charge $c$ and $G = \langle g \rangle$ a finite, cyclic group of automorphisms of $V$ of or\-der $n$. Then
\begin{enumerate}[i)]
\item
Let $q = e^{2 \pi i \tau}$. Then the twisted trace function $T(v,i,j,\tau)$ converges to a holomorphic function on $H$.
\item
For homogeneous $v\in V_{[k]}$ the twisted trace functions satisfy
\[  T(v,i,j,\gamma \tau)= \sigma(i,j,\gamma) (c\tau+d)^k T(v,(i,j)\gamma,\tau)  \]
for all $\gamma = \left( \begin{smallmatrix} a & b \\ c & d \end{smallmatrix} \right) \in \SL_2(\Z)$. The constants $\sigma(i,j,\gamma) \in \C$ depend only on $i$, $j$, $\gamma$ and the choice of the functions $\phi_{g^i}$.
\end{enumerate}
\end{thm}

\section{Simple currents} \label{simpcur}

In this section we show that a rational vertex operator algebra satisfying certain regularity conditions and whose modules are all simple currents has group-like fusion and that the conformal weights define a quadratic form on the fusion group. We also show that Zhu's representation is up to a character the Weil representation in this case.

\vspace*{2mm}
Let $V$ be a rational, $C_2$-cofinite vertex operator algebra of CFT-type. It is well-known that the fusion algebra $\V(V)$ of $V$ is a finite-dimensional, associative, commutative algebra over $\C$ with unit $V$. A $V$-module $M$ is called a simple current if the fusion product $M \boxtimes_V W$ is irreducible for  any irreducible $V$-module $W$. We define $\Irr(V) = \{ W^{\alpha} \, | \, \alpha \in D \}$ as the set of isomorphism classes of irreducible $V$-modules and assume that all $W^{\alpha}$ are simple currents. Then we can define a composition on $D$ by $W^{\alpha} \boxtimes_V W^{\beta} \cong W^{\alpha + \beta}$. The following result is well-known (see. \cite{LY}, Corollary 1).

\begin{thm} \label{thm.3.1}
Let $V$ be a simple, rational, $C_2$-cofinite, self-contragredient vertex operator algebra of CFT-type. Assume that all irreducible $V$-modules are simple currents. Then the fusion algebra $\V(V)$ of $V$ is the group algebra $\C[D]$ of some finite abelian group $D$, i.e.
\[
W^\alpha\boxtimes_V W^\beta\cong W^{\alpha+\beta}
\]
for all $\alpha,\beta\in D$. The neutral element is given by $W^0 = V$ and the inverse of $\alpha$ by $W^{-\alpha} = W^{\alpha'} \cong (W^{\alpha})'$.
\end{thm}

If $V$ is as in the theorem we say that $V$ has group-like fusion. We will assume this from now on. We denote the conformal weight of $W^{\alpha}$ by $\rho(W^{\alpha})$ and define
\[  q_{\Delta} : D \to \Q/\Z   \]
with
\[  q_{\Delta}(\alpha) = \rho( W^{\alpha} ) \! \mod 1    \]
and
\[   b_{\Delta} : D \times D \to \Q/\Z  \]
by 
\[  b_{\Delta} (\alpha,\beta) 
    = q_{\Delta}(\alpha + \beta) - q_{\Delta}(\alpha) - q_{\Delta}(\beta)  \! \mod 1 \, . \]

\begin{prp} \label{Smatcoeff}
Suppose $V$ has group-like fusion with fusion group $D$ and the modules $W^{\alpha}$ have positive conformal weights for $\alpha \neq 0$. Then 
\[  \mathcal{S}_{0 0} =  \mathcal{S}_{0 \alpha} = \frac{1}{\sqrt{|D|}}  \]
for all $\alpha \in D$.
\end{prp}
{\em Proof:}
Using the results on quantum dimensions in \cite{DJX} it is easy to see that $\mathcal{S}_{0 0} =  \mathcal{S}_{0 \alpha} = \mathcal{S}_{\alpha 0} \in \R^+$. Now
\[  1 = \delta_{0 0} = (\mathcal{S}^2)_{0 0} = 
\sum_{\gamma \in D} \mathcal{S}_{0 \gamma} \mathcal{S}_{\gamma 0} = |D|  \mathcal{S}_{0 0}^2  \]
implies the statement. \eop
 
\begin{prp}  \label{ops}
Suppose $V$ has group-like fusion with fusion group $D$ and the modules $W^{\alpha}$ have positive conformal weights for $\alpha \neq 0$. Then 
\[  \mathcal{S}_{\alpha \beta} = \mathcal{S}_{0 0} \, e( - b_{\Delta} (\alpha,\beta) )   \]
for all $\alpha, \beta \in D$.
\end{prp}
{\em Proof:}
The relation $TSTST=S$ in ${\SL}_2(\Z)$ implies 
$\mathcal{T}\mathcal{S}\mathcal{T}\mathcal{S}\mathcal{T} = \mathcal{S}$. 
Using the above formula for $\mathcal{T}$ we obtain
\[  \mathcal{S}_{\alpha \beta} = \sum_{\gamma \in D} \mathcal{S}_{\alpha \gamma} \mathcal{S}_{\gamma \beta} \,
                               e( q_{\Delta}(\alpha) + q_{\Delta}(\beta) + q_{\Delta}(\gamma) - c/8 )  \, . \]
By the Verlinde formula we have
\begin{align*} 
\mathcal{S}_{\alpha + \beta, \gamma} 
&= 
\sum_{\delta \in D} \mathcal{S}_{\delta \gamma} \delta_{\alpha + \beta, \delta}  
= \sum_{\delta \in D} \mathcal{S}_{\delta \gamma} \mathcal{N}_{\alpha \beta}^{\delta}  
= \sum_{\delta, \rho \in D} \mathcal{S}_{\delta \gamma} \frac{\mathcal{S}_{\alpha \rho} \mathcal{S}_{\beta \rho} \mathcal{S}_{-\delta, \rho}}{\mathcal{S}_{0 \rho}}  \\
&=
\sum_{\rho \in D} \frac{\mathcal{S}_{\alpha \rho} \mathcal{S}_{\beta \rho}}{\mathcal{S}_{0 \rho}} 
\sum_{\delta \in D} \mathcal{S}_{\delta \gamma} \mathcal{S}_{-\delta, \rho} 
= \sum_{\rho \in D} \frac{\mathcal{S}_{\alpha \rho} \mathcal{S}_{\beta \rho}}{\mathcal{S}_{0 \rho}} 
\, \delta_{\gamma \rho} 
=
\frac{\mathcal{S}_{\alpha \gamma} \mathcal{S}_{\beta \gamma}}{\mathcal{S}_{0 \gamma}} 
\end{align*}
because
\[  \sum_{\delta \in D} \mathcal{S}_{\delta \gamma} \mathcal{S}_{-\delta, \rho} =
\sum_{\delta, \mu \in D} \mathcal{S}_{\gamma \delta} \delta_{-\delta, \mu} \mathcal{S}_{\mu \rho} =
\sum_{\delta, \mu \in D} \mathcal{S}_{\gamma \delta} (\mathcal{S}^2)_{\delta \mu} \mathcal{S}_{\mu \rho} =
(\mathcal{S}^4)_{\gamma \rho} =
\delta_{\gamma \rho}  \, . \]
Hence
\[  \mathcal{S}_{\alpha + \beta, \gamma} \mathcal{S}_{0 \gamma}
= \mathcal{S}_{\alpha \gamma} \mathcal{S}_{\beta \gamma}  \]
so that
\begin{align*}
\mathcal{S}_{\alpha \beta} 
&= 
\sum_{\gamma \in D} \mathcal{S}_{\alpha \gamma} \mathcal{S}_{\gamma \beta} \,
                 e( q_{\Delta}(\alpha) + q_{\Delta}(\beta) + q_{\Delta}(\gamma) - c/8 ) \\
&= 
\sum_{\gamma \in D} \mathcal{S}_{\alpha + \beta, \gamma} \mathcal{S}_{0 \gamma} \, 
                 e( q_{\Delta}(\alpha) + q_{\Delta}(\beta) + q_{\Delta}(\gamma) - c/8 ) \\
&=
e( q_{\Delta}(\alpha) + q_{\Delta}(\beta) - c/8 )
\sum_{\gamma \in D} \mathcal{S}_{\alpha + \beta, \gamma} \mathcal{S}_{0 \gamma} \, e( q_{\Delta}(\gamma) ) \\
&= 
e( q_{\Delta}(\alpha) + q_{\Delta}(\beta) - c/8 )
\mathcal{S}_{\alpha + \beta, 0} \, e( -q_{\Delta}(\alpha + \beta) + c/8 ) \\[1mm]
&=
e( - b_{\Delta} (\alpha,\beta) ) \, \mathcal{S}_{\alpha + \beta, 0} 
\end{align*}
since
\[ 
\mathcal{S}_{\alpha 0} 
= 
\sum_{\gamma \in D} \mathcal{S}_{\alpha \gamma} \mathcal{S}_{\gamma 0} \,
                 e( q_{\Delta}(\alpha) + q_{\Delta}(\gamma) - c/8 )   
=
e( q_{\Delta}(\alpha) - c/8 )
\sum_{\gamma \in D} \mathcal{S}_{\alpha \gamma} \mathcal{S}_{\gamma 0} \,
                 e( q_{\Delta}(\gamma) ) \, . 
\]
This finishes the proof. \eop

\begin{thm} \label{thm.3.4}
Let $V$ be a simple, rational, $C_2$-cofinite, self-contragredient vertex operator algebra of CFT-type. Suppose $V$ has group-like fusion with fusion group $D$ and the modules $W^{\alpha}$ have positive conformal weights for $\alpha \neq 0$. Then 
\begin{align*} 
\mathcal{S}_{\alpha \beta} &= \frac{1}{\sqrt{|D|}} \, e( - b_{\Delta} (\alpha,\beta) ) , \\
\mathcal{T}_{\alpha \beta} &= e( q_{\Delta}(\alpha) -c/24 ) \delta_{\alpha \beta} .
\end{align*}
Moreover $q_{\Delta}$ is a quadratic form on $D$ and $b_{\Delta}$ the associated bilinear form.
\end{thm}
{\em Proof:} 
The formula for $\mathcal{S}$ follows from the previous two results. The relation  
$\mathcal{S}_{\alpha \gamma} \mathcal{S}_{\beta \gamma} =
\mathcal{S}_{\alpha + \beta, \gamma} \mathcal{S}_{0 \gamma}$ together with Proposition \ref{ops} shows that $b_{\Delta}$ is bilinear. We have $q_{\Delta}(0) = 0 \mod 1$ 
and $q_{\Delta}(\alpha) = q_{\Delta}(-\alpha)$ for all $\alpha \in D$. This implies that $q_{\Delta}$ is a quadratic form with associated bilinear form $b_{\Delta}$. \eop\\[-2mm]

\begin{prp}
Suppose $V$ has group-like fusion with fusion-group $D$ and the modules $W^{\alpha}$ have positive conformal weights for $\alpha \neq 0$. Then the bilinear form $b_{\Delta}$ is non-degenerate.
\end{prp} 
{\em Proof:} We have 
\[  \delta_{\alpha, -\beta} =  (\mathcal{S}^2)_{\alpha \beta} 
= \sum_{\gamma \in D} \mathcal{S}_{\alpha \gamma} \mathcal{S}_{\beta \gamma}
= \frac{1}{|D|} \sum_{\gamma \in D} e(-b_{\Delta} (\alpha + \beta, \gamma) ) \]
so that 
\[  \delta_{\alpha 0} = \frac{1}{|D|} \sum_{\gamma \in D} e(-b_{\Delta} (\alpha, \gamma) )  \, . \] 
This implies the statement. \eop\\[-2mm]

It follows that $D$ is a finite abelian group with a quadratic form $q_{\Delta} : D \to \Q/\Z$ such that the corresponding bilinear form is non-degenerate, i.e. a discriminant form \cite{N}. The group algebra $\C[D]$ carries two representations of the metaplectic group $\text{Mp}_2(\Z)$, one coming from Zhu's theorem and the other being the Weil representation  \cite{B3}, \cite{NRS}. These representations commute and the explicit formulas for $S$ and $T$ show that they are related by a character. As a consequence we have 

\begin{thm}
Let $V$ be a simple, rational, $C_2$-cofinite, self-contragredient vertex operator algebra of CFT-type. Suppose $V$ has group-like fusion with fusion group $D$ and the modules $W^{\alpha}$ have positive conformal weights for $\alpha \neq 0$. Then the central charge $c$ of $V$ is an integer and $D$ is a discriminant form under $q_{\Delta}$ of signature
\[  \sign(D) = c \mod 8 . \]
\end{thm}

\section{Abelian intertwining algebras}

In this section we show that the irreducible modules of a rational vertex operator algebra $V$ satisfying certain regularity conditions and with group-like fusion form an abelian intertwining algebra whose quadratic form is determined by the conformal weights. Restricting the sum to an isotropic subgroup gives a vertex operator algebra extending $V$.

\vspace*{2mm}
Let $D$ be a finite abelian group and $(F,\Omega)$ a normalised abelian $3$-cocycle on $D$ with coefficients in $\C^*$, i.e. the maps 
\begin{gather*}
F : D \times D \times D \to \C^* \\
\Omega : D \times D \to \C^* 
\end{gather*}
satisfy
\begin{gather*}
F(\alpha, \beta, \gamma) F(\alpha, \beta, \gamma + \delta)^{-1} F(\alpha, \beta + \gamma, \delta) F(\alpha + \beta, \gamma, \delta)^{-1} F(\beta, \gamma, \delta) = 1 \\
F(\alpha, \beta, \gamma)^{-1} \Omega(\alpha, \beta + \gamma) F(\beta, \gamma, \alpha)^{-1} 
= \Omega(\alpha, \beta) F(\beta, \alpha, \gamma)^{-1} \Omega(\alpha, \gamma) \\
F(\alpha, \beta, \gamma) \Omega(\alpha + \beta, \gamma) F(\gamma, \alpha, \beta)
= \Omega(\beta, \gamma) F(\alpha, \gamma, \beta) \Omega(\alpha, \gamma)
\end{gather*}
and
\begin{gather*} 
F(\alpha, \beta, 0) = F(\alpha, 0, \gamma) = F(0, \beta, \gamma) = 1 \\
\Omega(\alpha, 0 ) = \Omega( 0 , \beta) = 1 
\end{gather*}
for all $\alpha, \beta, \gamma, \delta \in D$. Define
\[ B : D \times D \times D \to \C^* \]
by
\[ B(\alpha, \beta, \gamma) = 
F(\beta, \alpha, \gamma)^{-1} \Omega(\alpha, \beta) F(\alpha, \beta, \gamma) \, . \]
We also define a quadratic form 
\[ q_{\Omega} : D \to \Q/\Z \]
by 
\[ \Omega(\alpha,\alpha) = e( q_{\Omega}(\alpha) )  \]
for all $\alpha \in D$. We denote the corresponding bilinear form by $b_{\Omega}$. The level of $(F,\Omega)$ is the smallest positive integer $N$ such that $N q_{\Omega}(\alpha) = 0 \! \mod 1$ for all $\alpha \in D$. 

An abelian intertwining algebra \cite{DL1} of level $N$ associated to $D, F$ and $\Omega$ is a $\C$-vector space $V$ with a $\frac{1}{N}\Z$-grading and a $D$-grading
\[ V = \bigoplus_{n \in \frac{1}{N}\Z} V_n = \bigoplus_{\alpha \in D} V^{\alpha}  \]
such that
\[  V^{\alpha} =  \bigoplus_{n \in \frac{1}{N}\Z} V_n^{\alpha} \, , \]
where $V_n^{\alpha} = V_n \cap V^{\alpha}$, equipped with a state-field correspondence
\begin{align*}
Y \, : \, V \, & \, \to \End (V) [[ z^{-1/N}, z^{1/N} ]]  \\
          a \, & \, \mapsto Y(a,z) = \sum_{n \in \frac{1}{N} \Z} a_n z^{-n-1} 
\end{align*}
and with two distinguished vectors $\vac \in V_0^0$, $ \omega \in V_2^0$ satisfying the following conditions. For $\alpha, \beta \in D$, $a,b \in V$ and $n \in \frac{1}{N}\Z$ 
\begin{gather*}
a_n V^{\beta} \subset V^{\alpha + \beta}  \quad \text{if $a \in V^{\alpha}$} \\ 
a_nb = 0  \quad \text{for $n$ sufficiently large} \\
Y(a,z) \vac \in V[[z]] \quad \text{and} \quad \lim_{z \to 0} Y(a,z) \vac = a \\
Y(a,z)|_{V^{\beta}} = \sum_{n \in b_{\Omega}(\alpha, \beta) + \Z} a_n z^{-n-1} 
\quad \text{if $a \in V^{\alpha}$,}
\end{gather*} 
the Jacobi identity
\begin{multline*}  
x^{-1} \left( \frac{y-z}{x} \right)^{b_{\Omega}(\alpha, \beta)} 
\delta\left( \frac{y-z}{x} \right)
Y(a,y) Y(b,z) c \\
- B(\alpha, \beta , \gamma) \, x^{-1} \left( \frac{z-y}{e^{\pi i}x} \right)^{b_{\Omega}(\alpha, \beta)} 
\delta\left( \frac{z-y}{-x} \right)
Y(b,z) Y(a,y) c \\
= F(\alpha, \beta , \gamma) z^{-1}  \delta\left( \frac{y-x}{z} \right) 
Y(Y(a,x)b,z) \left( \frac{y-x}{z} \right)^{-b_{\Omega}(\alpha, \gamma)}c
\end{multline*} 
holds for all $a \in V^{\alpha}$,  $ b \in V^{\beta}$ and $c \in V^{\gamma}$, 
the operators $L_n$ defined by 
\[  Y(\omega,z) = \sum_{n \in \Z} L_n z^{-n-2}  \]
satisfy
\[  [L_m,L_n] = (m-n) L_{m+n} + \frac{m^3-m}{12} \delta_{m+n} c  \]
for some $c \in \C$ and 
\begin{gather*}
L_0a = na \quad \text{for $a \in V_n$,}  \\
\frac{d}{dz} Y(a,z) = Y(L_{-1}a,z) \, .    
\end{gather*} 

The cohomology class of $(F,\Omega)$ in $H^3_{\text{ab}}(G,\C^*)$ is determined by $q_{\Omega}$ (see.\ \cite{ML}, Theorem 3 and \cite{DL1}, Remark 12.22). We also remark that rescaling the intertwining operators amounts to changing $(F,\Omega)$ by a coboundary (\cite{DL1}, Remark 12.23).

A consequence of the Jacobi identity is the skew-symmetry formula
\[   Y(a,z)b = \frac{1}{\Omega(\bt,\al)} \, e^{zL_{-1}} Y(b,e^{-\pi i} z)a    \]
for $a \in V^{\al}$, $b \in V^{\bt}$.

It is well-known that the vertex operator algebra of a positive definite, even lattice together with its irreducible modules forms an abelian intertwining algebra where $q_{\Omega}$ is determined by the conformal weights (\cite{DL1}, Theorem 12.24). More generally we have 

\begin{thm}\label{thm.4.7}
Let $V$ be a simple, rational, $C_2$-cofinite, self-contragredient vertex operator algebra of CFT-type. Assume that $V$ has group-like fusion, i.e. the fusion algebra of $V$ is $\C[D]$ for some finite abelian group $D$, and that the irreducible $V$-modules $W^{\alpha}$, $\alpha \neq 0$ have positive conformal weights. Then 
\[  W = \bigoplus_{\alpha \in D} W^{\alpha}  \]
can be given the structure of an abelian intertwining algebra with normalised abelian $3$-cocycle $(F,\Omega)$ such that
\[ q_{\Omega} = - q_{\Delta} \, . \]
\end{thm}
{\em Proof:} 
Choosing non-trivial intertwining operators between the modules $W^{\alpha}$ of $V$ the Jacobi identity defines maps $F$ and $\Omega$. Huang has shown that $(F,\Omega)$ is an abelian $3$-cocycle on $D$ (see.\ \cite{H2}, Theorem 3.7 and \cite{H1}). Any abelian $3$-cocycle is cohomologous to a normalised one so that $W$ is an abelian intertwining algebra. Furthermore the modules of $V$ form a modular tensor category with the twist morphisms $\theta_{\alpha} : W^{\alpha} \to W^{\alpha}$ on the irreducible modules given by $\theta_{\alpha} = e(q_{\Delta}(\alpha))\id_{W^\alpha}$ (see. \cite{H4}). Since $V$ has group-like fusion the braiding iso\-morphism 
\[ c_{\alpha,\beta}: W^{\al} \boxtimes W^{\bt} \to W^{\bt} \boxtimes W^{\al}  \]
is for $\alpha = \beta$ given by
\[  c_{\alpha,\alpha} = 
     e(q_\Omega(\alpha)) e(2q_\Delta(\alpha)) \id_{W^{\alpha} \boxtimes W^{\alpha}}  \, . \]
We have
\[ \tr_{W^{\alpha}} \theta_{\alpha} = 
         \tr_{W^{\alpha} \boxtimes W^{\alpha}} c_{ \alpha, \alpha}  \]
(see.\ Proposition 2.32 in \cite{DGNO}) so that
\[  e(q_{\Delta}(\alpha)) \tr_{W^{\alpha}} \id_{W^{\alpha}} = 
e(q_\Omega(\alpha)) e(2q_\Delta(\alpha)) \tr_{W^{2\alpha}} \id_{W^{2\alpha}}    \, . \]
The trace 
\[ d_{\alpha} = \tr_{W^\alpha} \id_{W^\alpha} \]
is the categorical quantum dimension of $W^{\alpha}$. For Huang's construction under the positivity assumption on the conformal weights it coincides with the definition of the quantum dimension as the limit of a certain character ratio 
so that $d_{\alpha} = 1$ if $W^{\alpha}$ is a simple current (see.\ \cite{DLN}, Proposition 3.11 and \cite{DJX}, Proposition 4.17). Then the above identity implies
\[  e(q_{\Omega}(\alpha)) = e( - q_{\Delta}(\alpha)) \]
for all $\al \in D$. This proves the theorem. \eop\\[-2mm]

This result is also stated as Theorem 2.7 in \cite{HS}. In the above proof we fill a gap pointed out by S.\ Carnahan. The theorem can also be derived from Bantay's formula for the Frobenius-Schur indicator and the fact that the theorem holds for abelian intertwining algebras associated with positive definite, even lattices. 

An explicit construction of the 3-cocycle $(F,\Omega)$ in the lattice case is described in \cite{DL1}, Theorem 12.24. Since every discriminant form can be realised as dual quotient of a positive definite, even lattice this also covers the general case.

As an application we obtain

\begin{thm} \label{extensionthm}
Let $V$ be a simple, rational, $C_2$-cofinite, self-contragredient vertex operator algebra of CFT-type with group-like fusion and fusion group $D$. Suppose that the irreducible $V$-modules $W^{\alpha}$, $\alpha \in D \setminus \{ 0 \}$ have positive conformal weights. Let $H$ be an isotropic subgroup of $D$ with respect to $q_{\Delta}$. Then  
\[  W^H = \bigoplus_{\gamma \in H} W^{\gamma}  \]
admits the structure of a simple, rational, $C_2$-cofinite, self-contragredient vertex operator algebra of CFT-type extending the vertex operator algebra structure on $V$. If $H=H^{\perp}$ then $W^H$ is holomorphic.
\end{thm}
{\em Proof:}
Suppose $q_{\Delta}|_H=q_{\Omega}|_H = 0 \! \mod 1$. Then $(F|_H,\Omega|_H)$ is cohomologous to the trivial $3$-cocycle in $H^3_{\text{ab}}(H,\C^*)$. Hence the abelian intertwining algebra $W^H$ admits the structure of a vertex operator algebra upon rescaling of the intertwining operators. The irreducible modules of $W^H$ are given by 
\[  W^{H, \gamma} = \bigoplus_{\alpha \in \gamma + H} W^{\alpha}  \] 
where $\gamma$ ranges over $H^{\perp}/H$ (see.\ \cite{Y}, Theorems 3.2 and 3.3). \eop

\section{Orbifolds}\label{Section.orb}

Let $G = \langle g \rangle$ be a finite, cyclic group of order $n$ acting on a holomorphic vertex operator algebra $V$. We show that the fixed-point subalgebra $V^G$ has group-like fusion and that the fusion group is a central extension of $\Z_n$ by $\Z_n$ whose isomorphism type is fixed by the conformal weight of the twisted module $V(g)$. We also determine the $\mathcal{S}$-matrix and describe the level of the trace functions. If the twisted modules $V(g^j)$ have positive conformal weights for $j \neq 0$ then the direct sum of the irreducible $V^G$-modules is an abelian intertwining algebra and the restriction of the sum to an isotropic subgroup with respect to the conformal weights is a vertex operator algebra extending $V^G$. Our approach is inspired by Miyamoto's theory of $\Z_3$-orbifolds \cite{M1}.

\vspace*{2mm}
Let $V$ be a simple, rational, $C_2$-cofinite, holomorphic vertex operator algebra of CFT-type and $G = \langle g \rangle$ a finite, cyclic group of order $n$ acting by auto\-mor\-phisms on $V$. Then we have (see.\ \cite{CM}, \cite{M2})

\begin{thm} \label{ratvg}
The fixed-point subalgebra $V^G$ is a simple, rational, $C_2$-cofinite, self-contra\-gre\-dient vertex operator algebra of CFT-type. Every irreducible $V^G$-module is isomorphic to a $V^G$-submodule of the irreducible $g^i$-twisted $V$-module $V(g^i)$ for some $i$.
\end{thm}

Recall that for each $h\in G$ there is a representation $\phi_h:G \rightarrow \Aut_{\C}(V(h))$ of $G$ on the vector space $V(h)$ such that $\phi_h(k)Y_{V(h)}(v,z)\phi_h(k)^{-1} = Y_{V(h)}(kv,z)$ for all $k\in G$ and $v\in V$ and these representations are unique up to multiplication by an $n$-th root of unity. If $h$ is the identity we can and will assume that $\phi_h(k) = k$ for all $k \in G$. We will often write $\phi_i$ for $\phi_{g^i}$. We denote the eigen\-space of $\phi_i(g)$ in $V(g^i)$ with eigenvalue $e(j/n)=e^{2 \pi i j/n}$ by $W^{(i,j)}$, i.e. 
\[ W^{(i,j)} =\{ v \in V(g^i) \, | \, \phi_i(g) v = e(j/n) v \} \, . \]
The twisted trace functions $T(v,i,j,\tau)$ transform under $S$ as
\[  T(v,i,j,S \tau) = T(v,i,j,-1/\tau) = \tau^{\wt[v]} \lambda_{i,j} T(v,j,-i,\tau) \]
where the $\lambda_{i,j} = \sigma(i,j,S)$ are complex numbers depending only on $i,j$.

Combining Theorem \ref{ratvg} with Theorem 2 of \cite{MT} we obtain

\begin{thm}
The fixed-point subalgebra $V^G$ has exactly $n^2$ distinct irreducible modules up to isomorphism, the eigenspaces $W^{(i,j)}$.
\end{thm}

We describe the contragredient module ${W^{(i,j)}}'$ of $W^{(i,j)}$.

\begin{prp}
We have
\[  {W^{(i,j)}}' \cong W^{(-i,\alpha(i)-j)}  \]
for some function $\al: \Z_n \to \Z_n$ satisfying $\al(i) = \al(-i)$ and $\al(0) = 0$. 
\end{prp}
{\em Proof:}
Let $W$ be an irreducible twisted or untwisted $V$-module and let $\phi_W : G \to \Aut_{\C}(W)$ denote the action of $G$ on the vector space $W$ such that 
\[   Y_W(g v,z) = \phi_W(g) Y_W(v,z) \phi_W(g)^{-1}    \]
for all $g \in G$. The definition of the contragredient module implies  
\[  Y_{W'}(g v,z) = \phi_W'(g)^{-1} Y_{W'}(v,z) \phi_W'(g)  \]
where $\phi_W': G \to \Aut_{\C}(W')$ denotes the graded dual representation of $\phi_W$. Hence by uniqueness the representation $\phi_{W'}$ is proportional to the inverse of $\phi_W'$. From $V(g^i)' \cong V(g^{-i})$ we obtain 
\[  \bigoplus_{j \in \Z_n}{W^{(i,j)}}' \cong \bigoplus_{j \in \Z_n} W^{(-i,j)}  \]
so that ${W^{(i,j)}}' \cong W^{(-i, j')}$ for some $j' \in \Z_n$. Write $\phi_{-i} = e(\al(i)/n) {\phi_i'}^{-1}$ for some $\al(i) \in \Z_n$. The subspace $W^{(i,j)}$ of $V(g^i)$ is the eigenspace of $\phi_i(g)$ with eigenvalue $e(j/n)$. The subspace ${W^{(i,j)}}' \cong W^{(-i, j')}$ of $V(g^i)' \cong V(g^{-i})$ is the eigenspace of ${\phi_i'}(g) = e(\al(i)/n) \phi_{-i}^{-1}(g)$ with the same eigenvalue $e(j/n)$. Hence $e(j/n) = e(\al(i)/n) e(-j'/n)$ so that
\[   {W^{(i,j)}}' \cong W^{(-i,\alpha(i)-j)} \, . \]
Now $W'' \cong W$ implies $j=\alpha(-i)-(\alpha(i)-j)$ and $\alpha(i)=\alpha(-i)$. With $\phi_0(g)=g$ we have $W^{(0,0)}=V^G$ so that $\alpha(0)=0$. \eop\\[-2mm]

We determine the fusion algebra $\V(V^G)$ of $V^G$.

\begin{prp} \label{Smatrix.to.lambda}
The $\mathcal{S}$-matrix of $V^G$ is given by
\[  \mathcal{S}_{(i,j),(k,l)} = \frac{1}{n} e(-(il+jk)/n) \lambda_{i,k}   \, . \]
\end{prp}
{\em Proof:}
We have 
\begin{gather*}
T(v,i,j,\tau) = \tr_{V(g^i)} o(v) \phi_i(g^j) q^{L_0-c/24}  \\[1mm]
= \sum_{k \in \Z_n} \tr_{W^{(i,k)}} o(v) e(jk/n) q^{L_0-c/24} 
= \sum_{k\in\Z_n} e(jk/n) T_{W^{(i,k)}} (v,\tau)
\end{gather*}
so that
\[  T_{W^{(i,j)}}(v,\tau) = \frac{1}{n} \sum_{l\in\Z_n} e(-lj/n) T(v,i,l,\tau)  \]
and
\begin{align*}
T_{W^{(i,j)}}(v,S\tau) 
&= \frac{1}{n} \sum_{k \in \Z_n} e(-jk/n) T(v,i,k,S\tau)  \\
&= \frac{1}{n} \sum_{k \in \Z_n} e(-jk/n) \tau^{\wt[v]} \lambda_{i,k} T(v,k,-i,\tau) \\ 
&= \frac{1}{n} \sum_{k,l \in \Z_n} e(-(jk+il)/n) \tau^{\wt[v]} \lambda_{i,k} T_{W^{(k,l)}}(v,\tau) 
\end{align*} 
by twisted modular invariance. \eop

\begin{prp} \label{proplambda}
The constants $\lambda_{i,j}$ satisfy
\begin{align*}
\lambda_{i,j}               & = \lambda_{j,i} \\
\lambda_{i,j} \lambda_{i,-j} & = e(i \al(j)/n) \\
\lambda_{0,i}               & = 1
\end{align*}
for all $i,j \in \Z_n$.
\end{prp}
{\em Proof:}
The first equation follows from the symmetry of the $\mathcal{S}$-matrix. $\mathcal{S}^2$ is a permutation matrix sending the index $(i,j)$ to the index of the contragredient module $(-i,-j+\al(i))$, i.e.
\[ (\mathcal{S}^2)_{(i,j),(l,m)} = \delta_{i,-l} \delta_{j, -m +\al(i)} \, . \]
By the previous proposition we have
\begin{align*} 
(\mathcal{S}^2)_{(i,j),(l,m)} 
& = \sum_{a,b \in \Z_n} \mathcal{S}_{(i,j),(a,b)} \mathcal{S}_{(a,b),(l,m)}  \\
& = \frac{1}{n^2} \sum_{a,b\in \Z_n} e(-(ib+ja+am+bl)/n) \lambda_{i,a} \lambda_{a,l} \\
& = \frac{1}{n^2} \sum_{a\in \Z_n} e(-a(j+m)/n) \lambda_{i,a} \lambda_{a,l} 
                  \sum_{b\in \Z_n} e(-b(i+l)/n) \\
& = \delta_{i,-l} \, \frac{1}{n} \sum_{a\in \Z_n} e(-a(j+m)/n) \lambda_{i,a} \lambda_{a,l} 
\end{align*}
so that 
\[ \frac{1}{n} \sum_{a\in \Z_n} e(-a(j+m)/n) \lambda_{i,a} \lambda_{a,-i} 
    = \delta_{j,\al(i)-m}  \]
and
\[ \frac{1}{n} \sum_{a\in \Z_n} e(-ab/n) \lambda_{i,a} \lambda_{a,-i} 
    = \delta_{b-m,\al(i)-m} \, . \]
Multiplying with $e(db/n)$ and summing over $b$ gives the second equation of the proposition. This equation implies $\lambda_{0,i}^2 = 1$. In order to prove the last equation we show that $\lambda_{0,i} \in \R_{\geq 0}$. We have
$T( \vac , 0, j,-1/\tau) = \lambda_{0,j} T( \vac ,j,0,\tau)$ so that
\[ \sum_{k=0}^\infty e((-1/\tau)(k-c/24)) \tr_{V_k} g^j  =  
\lambda_{0,j} \sum_{ k \in \rho_j+(1/n)\Z_{\geq 0}} e((k-c/24)\tau) \dim( V(g^j)_k ) 
\]
where $\rho_j \in \Q$ is the conformal weight of $V(g^j)$ and $c \in \Q$ the central charge of $V$. Specialising to $\tau = it$ with $t \in \R_{>0}$ we obtain
\[  \frac{1}{\lambda_{0,j}} \sum_{k=0}^{\infty} e^{-2 \pi k/t}\tr_{V_k} g^j = 
e^{2 \pi c (t-1/t)/24} \sum_{k \in \rho_j+(1/n)\Z_{\geq 0}} e^{-2\pi k t} \dim(V(g^j)_k) 
\in \R_{\geq 0} \, . \]
The limit of the left hand side for $t \to 0$ exists and is $1/\lambda_{0,j}$ because $V$ is of CFT-type. Hence $\lambda_{0,j}$ is a non-negative real number. \eop

\begin{prp} \label{prop.mods.are.simple}
The irreducible $V^G$-modules $W^{(i,j)}$ are simple currents. 
\end{prp}
{\em Proof:}
As before let $(i,j)'$ denote the index of the contragredient module of $W^{(i,j)}$. Then
\begin{align*}
\mathcal{S}_{(i,j),(k,l)} \mathcal{S}_{(i,j)',(k,l)} 
&= \mathcal{S}_{(i,j),(k,l)} \mathcal{S}_{(-i,\alpha(i)-j),(k,l)}   \\
&= \frac{1}{n} e(-(kj+il)/n) \lambda_{i,k} \frac{1}{n} e(-(k(\alpha(i)-j)+(-i)l)/n) \lambda_{-i,k} \\
&= \frac{1}{n^2} e(-k\alpha(i)/n) \lambda_{i,k} \lambda_{-i,k} = \frac{1}{n^2}
\end{align*}
by Proposition \ref{proplambda}. We compute the fusion coefficients with the Verlinde formula
\begin{align*}
N_{(i,j),(i,j)'}^{(l,k)}
&= \sum_{a,b \in \Z_n} \frac{\S_{(i,j),(a,b)} \S_{(i,j)',(a,b)} \S_{(a,b),(l,k)'} }{ \S_{(0,0),(a,b)} }  \\
&= \frac{1}{n^2} \sum_{a,b \in \Z_n} \frac{ \S_{(a,b),(l,k)'} }{ \S_{(0,0),(a,b)} }  
= \frac{1}{n^2} \sum_{a,b \in \Z_n} \frac{ \S_{(a,b),(-l,\alpha(l)-k)}}{\S_{(0,0),(a,b)}}  \\[1mm]
&= \frac{1}{n^2} \sum_{a,b \in \Z_n} 
   \frac{ e(-((-l)b+a(\alpha(l)-k))/n) \lambda_{a,-l} }{ \lambda_{0,a} }  \\
&= \frac{1}{n^2} \sum_{a \in \Z_n} e(a(k-\alpha(l))/n) \lambda_{a,-l} \sum_{b \in \Z_n} e(lb/n)  \\[1mm]
&= \delta_{l,0}\frac{1}{n} \sum_{a \in \Z_n} e(a(k-\alpha(0))/n) \lambda_{a,0} 
= \delta_{l,0} \frac{1}{n} \sum_{a\in\Z_n} e(ak/n) 
= \delta_{l,0}\delta_{k,0} \, .
\end{align*}
This means 
\[  W^{(i,j)} \boxtimes W^{(i,j)'} \cong W^{(0,0)}  \]
for all $i,j \in \Z_n$. By Corollary 1 in \cite{LY} this implies that all $W^{(i,j)}$ are simple currents. \eop\\[-2mm]

Hence the fusion algebra of $V^G$ is the group algebra $\C[D]$ of a finite abelian group $D$ of order $n^2$. Propositions \ref{Smatrix.to.lambda} and \ref{proplambda} show 
\[  \mathcal{S}_{(0,0),(0,0)} = \mathcal{S}_{(0,0),(i,j)} = \frac{1}{\sqrt{|D|}}  \, . \]
Since the positivity assumption in the results of Section \ref{simpcur} only enters through Proposition \ref{Smatcoeff} this implies that all results of Section \ref{simpcur} hold for $V^G$ without the assumption on the conformal weights. In particular the reduction of the conformal weights modulo $1$ defines a quadratic form $q_{\Delta}$ on $D$ whose associated bilinear form $b_{\Delta}$ is non-degenerate.  

\begin{prp} \label{teahupoo}
The fusion product takes the form  
\[  W^{(i,j)} \boxtimes W^{(k,l)} \cong W^{(i+k,j+l+c(i,k))}   \]
for some symmetric, normalised 2-cocycle $c: \Z_n \times \Z_n \to \Z_n$ satisfying
\[  e(-ac(i,k)/n) = \frac{ \lambda_{i,a}\lambda_{k,a} }{ \lambda_{i+k,a} }   \]
for all $i,k,a \in \Z_n$.
\end{prp} 
{\em Proof:}
Let $W^{(i,j)} \boxtimes W^{(k,l)} \cong W^{(s,t)}$. Then 
\[  \mathcal{S}_{(i,j),(a,b)} \mathcal{S}_{(k,l),(a,b)} = \frac{1}{n} \mathcal{S}_{(s,t),(a,b)}   \]
for all $a,b \in \Z_n$ as shown in the proof of Proposition \ref{ops}. This implies
\[  \lambda_{i,a} \lambda_{k,a}/ \lambda_{s,a} = e(-(sb+ta)/n) e((ib+ja)/n) e((kb+la)/n)  \, .  \]
Taking $a=0$ we obtain $s = i+k \! \mod n$ because $\lambda_{i,0} =1$. For $b=0$ we get 
\[  \lambda_{i,a} \lambda_{k,a}/ \lambda_{i+k,a} = e((j+l-t)a/n)  \, . \]
This shows that $t-j-l$ depends only on $i$ and $k$ and we define $c(i,k) = t-j-l$. The associativity of the fusion algebra $\V(V^G)$ implies that $c: \Z_n\times \Z_n \to \Z_n$ is a $2$-cocycle.
This cocycle is symmetric since $\V(V^G)$ is commutative and normalised because $\lambda_{0,1}=1$. \eop \\

The 2-cocycle $c$ is related to the map $\alpha$ by
\[  \alpha(i) + c(i,-i) = 0 \mod n  \, . \]

The maps $\Z_n \to D, \, j \mapsto (0,j)$ and $D \to \Z_n, \, (i,j) \mapsto i$  give an exact sequence
\[  0 \to \Z_n \to D \to \Z_n \to 0 \, ,   \]
i.e.\ $D$ is a central extension of $\Z_n$ by $\Z_n$. This extension is determined up to isomorphism by the cohomology class of the $2$-cocycle $c$ in $H^2(\Z_n,\Z_n)$. The group $H^2(\Z_n,\Z_n)$ is isomorphic to $\Z_n$ and the $2$-cocycle $c_d$ corresponding to $d$ in $\Z_n$ is represented by 
\[  c_d(i,j) = 
\begin{cases} \, 0 & \text{if } i_n + j_n < n    \\  
              \, d & \text{if } i_n + j_n \geq n \end{cases}  \]
where $i_n$ denotes the representative of $i$ in $\{ 0, \ldots, n-1 \}$. Since $c$ is normalised its cohomology class $d$ can be determined by
\[  d = c(1,1) + c(1,2) + \ldots + c(1,n-1)  \mod n \, . \]

We write $\rho_i$ for the conformal weight $\rho(V(g^i))$ of the irreducible $g^i$-twisted $V$-module $V(g^i)$.

\begin{prp}  \label{choicephi}
We can choose the representation $\phi_i$ of $G$ on $V(g^i)$ such that
\[   \phi_i(g^{(i,n)}) = e\big( (i/(i,n))^{-1} ( L_0 - \rho_i)  \big)   \]
where $(i/(i,n))^{-1}$ denotes the inverse of $i/(i,n)$ modulo $n/(i,n)$. 
\end{prp}
{\em Proof:} First we consider the case $(i,n)=1$. The module $W = V(g^i)$ has $L_0$-decomposition 
\[  W = \bigoplus_{j \in \Z} W_{\rho_i + j/n} \, . \]
For $k \in \Z_n$ define 
\[  W_k = \bigoplus_{\substack{j \in \Z \\[0.4mm] j = k \bmod n}} W_{\rho_i + j/n} \, . \]
Let $v \in V$ be a homogeneous eigenvector of $g$ with eigenvalue $e(r/n)$, i.e. $gv = e(r/n)v$. Then by the definition of $g^i$-twisted modules the modes $v_t$ of $Y_W(v,z)$ vanish unless $t \in -ir/n + \Z$ and increase the $L_0$-eigenvalue by $\wt(v)-(t+1) \in ir/n + \Z$. Hence
$v_t W_k \subset W_{k + ir}$. Define the map 
\[   \xi(g^j) = e( i^{-1}j L_0 )   \]
where $i^{-1}$ denotes the inverse of $i$ modulo $n$. Then for homogeneous $w \in W$ and $t \in -ir/n + \Z$ we have $(g^jv)_tw = e(jr/n) v_tw$ and
\begin{align*}
\xi(g^j) v_t \xi(g^j)^{-1} w 
&= e( i^{-1}j L_0 ) v_t e( - i^{-1}j L_0 ) w 
= e( i^{-1}j L_0 ) v_tw e( - i^{-1}j \wt(w) ) \\
&= e( i^{-1} j(\wt(w) + ir/n) ) v_tw e( - i^{-1}j \wt(w) ) 
= e( jr/n)v_tw  
\end{align*}
so that
\[  \xi(g^j) Y_W (v,z) \xi(g^j)^{-1} w = Y_W(g^jv,z) w  \, . \]
The same relation holds for
\[   \phi(g^j) = e( i^{-1}j (L_0 - \rho_i) )   \, . \]
The shift by $\rho_i$ implies that $\phi$ defines a group homomorphism $\phi : G \rightarrow \Aut_{\C}(W)$. Hence by uniqueness $\phi$ is up to a constant equal to $\phi_i$ or put differently we can choose $\phi_i$ as $\phi$. This finishes the proof for $(i,n)=1$. The general case follows by replacing $g$ with $g^{(i,n)}$. \eop

\begin{prp}  
Choose the representation $\phi_i$ of $G$ on $V(g^i)$ such that 
\[   \phi_i(g^{(i,n)}) = e( (i/(i,n))^{-1} ( L_0 - \rho_i) ) \, . \]
Then the irreducible $V^G$-modules $W^{(i,j)}$ have conformal weights 
\[   \rho(W^{(i,j)}) = \rho_i + ij/n  \mod 1    \]
and the function $\alpha$ satisfies
\[   i \alpha(i) = 0  \mod n  \, . \]
\end{prp}
{\em Proof:} The first statement is easy to see. For the second statement about the map $\al$ observe that $\rho_i = \rho_{-i}$. \eop\\[-2mm]

If we choose the maps $\phi_i$ as in the proposition then the quadratic form $q_{\Delta}$ on $D$ is given by 
\[   q_{\Delta}((i,j)) = \rho_i + ij/n \mod 1 \]
and
\[   b_{\Delta}((i,j),(k,l)) 
= \rho_{i+k} - \rho_i - \rho_k + \frac{il+jk}{n} + \frac{(i+k)c(i,k)}{n}  \mod 1  \, . \]
We will use this choice in the proofs of the next $3$ results.

\begin{prp} \label{iggypop}
The cohomology class of the central extension defined by the $2$-cocycle $c$ is given by
\[  d = 2 n^2 \rho_1  \mod n  \, . \]
\end{prp}
{\em Proof:}
We have $q_{\Delta}((i,j)) = \rho_i + ij/n \! \mod 1$ so that
\[  \lambda_{i,i} 
= n e(2ij/n) \mathcal{S}_{(i,j),(i,j)}
= e(2ij/n) e(- 2q_{\Delta}((i,j)))
= e(-2 \rho_i )  \]
by Proposition \ref{Smatrix.to.lambda} and Theorem \ref{thm.3.4}. It follows that
\begin{align*}
e(d/n) 
& = e((c(1,1)+c(1,2)+\ldots+c(1,n-1))/n) \\
& = \frac{\lambda_{2,1}}{\lambda_{1,1}\lambda_{1,1}} \, \frac{\lambda_{3,1}}{\lambda_{1,1}\lambda_{2,1}} 
     \cdot \ldots \cdot \frac{\lambda_{n,1}}{\lambda_{1,1}\lambda_{n-1,1}} 
= \frac{\lambda_{n,1}}{\lambda_{1,1}^n} = \lambda_{1,1}^{-n} = e( 2 n \rho_1 )
\end{align*}
because $\lambda_{n,1} = \lambda_{0,1} = 1$.  \eop\\[-2mm]

Since $D$ has order $n^2$ the bilinear form $b_{\Delta}$ takes its values in $(1/n^2)\Z/\Z$ and the associated quadratic form $q_{\Delta}$ in $(1/2n^2)\Z/\Z$. We show that the values of $q_{\Delta}$ actually lie in $(1/n^2)\Z/\Z$.

\begin{thm} \label{rho1.values}
The unique irreducible $g$-twisted $V$-module $V(g)$ has conformal weight 
\[  \rho_1  \in (1/n^2)\Z   \]
and more generally $V(g^i)$ has conformal weight $\rho_i \in ((n,i)^2/n^2)\Z$. 
\end{thm}
{\em Proof:} 
It is sufficient to prove the statement for $i=1$. We have
\[ n^2 \rho_1 = n^2 q_{\Delta}((1,j)) =  q_{\Delta}(n(1,j)) = q_{\Delta}((0,k))  \mod 1 \]
for some $k \in \Z_n$. But this last value is $0 \! \mod 1$. This proves the theorem. \eop\\[-2mm]

This result generalises Theorem 1.6 (i) in \cite{DLM}. The value of $\rho_1$ determines the group structure of $D$. We will see that it also determines the quadratic form $q_{\Delta}$ up to isomorphism.

We define the type $t \in \Z_n$ of $g$ by
\[  t = n^2 \rho_1 \mod n  \, . \]
Then $d = 2t \! \mod n$. Let $N$ be the smallest positive multiple of $n$ such that $N \rho_1 = 0 \! \mod 1$, i.e.\ $N = n^2/(t,n)$. Then $N$ is the level of $q_{\Delta}$. 

\begin{prp}
The conformal weights $\rho_i$ satisfy 
\[ \rho_i = \frac{i^2 t}{n^2} \mod \, \frac{(i,n)}{n}  \, . \]
\end{prp}
{\em Proof:}
We have 
\begin{align*}
i^2 \rho_1 
&= i^2 q_{\Delta}((1,0)) =  q_{\Delta}(i(1,0))  
= q_{\Delta}( ( i , c(1,1) + \ldots + c(1,i-1) ) \\[1mm]
&= \rho_i + \frac{i(c(1,1) + \ldots + c(1,i-1))}{n}  
\end{align*}
so that 
\[  i^2 \rho_1 = \rho_i \mod \, \frac{(i,n)}{n}   \, . \]
On the other hand $i^2 \rho_1 = i^2t/n^2  \! \mod (i,n)/n$. This proves the proposition. \eop\\[-2mm]

The $2$-cocycle $c$ is cohomologous to $c_d$ with $d = 2t \! \mod n$. Changing the representations $\phi_i$ amounts to shifting $c$ by a coboundary. 

Recall that $i_n$ denotes the representative of $i \in \Z_n$ in $\{ 0, \ldots, n-1 \}$.
\begin{thm} \label{thm.5.12}
Let $V$ be a simple, rational, $C_2$-cofinite, holomorphic vertex operator algebra of CFT-type and $G = \langle g \rangle$ a finite, cyclic group of automorphisms of $V$ of order $n$. Suppose $g$ has type $t \! \mod n$ and let $d = 2t \! \mod n$. Then we can define the maps $\phi_i$ such that
\begin{enumerate}[i)]
 \item $W^{(i,j)} \boxtimes W^{(k,l)} \cong W^{(i+k,j+l+c_{d}(i,k))}$
 \item $W^{(i,j)}$ has conformal weight $q_{\Delta}((i,j)) = \dfrac{ij}{n} + \dfrac{i_n^2t_n}{n^2} \mod 1$
 \item ${W^{(i,j)}}'\cong W^{(-i,-j-c_{d}(i,-i))}$
 \item $\mathcal{S}_{(i,j),(k,l)} = \dfrac{1}{n} e(-(il+jk)/n) \lambda_{i,k} 
                      = \dfrac{1}{n} e(-(il+jk)/n) e(-2t_ni_nk_n/n^2)$, \\[2mm]
       i.e.\ $\lambda_{i,k} = e(-2t_n i_n k_n/n^2)$
\end{enumerate}
for $i,j,k,l \in \Z_n$.
\end{thm}
{\em Proof:} We prove the statement for $t = 0 \! \mod n$. The general case is similar. Then $\rho_i = 0 \! \mod (i,n)/n$. This implies that we can choose the representations $\phi_i$ of $G$ on $V(g^i)$ such that
\[   \phi_i(g^{(i,n)}) = e\big( (i/(i,n))^{-1} L_0 \big)   \]
(see.\ the proof of Proposition \ref{choicephi}). Then
\[  q_{\Delta}((i,j)) = \rho(W^{(i,j)}) = \frac{ij}{n} \mod 1 \]
and
\[  \lambda_{i,i} = n e(2ij/n) \mathcal{S}_{(i,j),(i,j)} 
                 = e(2ij/n) e(- 2q_{\Delta}((i,j))) = 1  \]
(see.\ the proof of Proposition \ref{iggypop}). The equation
\[  \frac{ \lambda_{i,a}\lambda_{k,a} }{ \lambda_{i+k,a} } = e(-ac(i,k)/n)  \]
in Proposition \ref{teahupoo} implies that the $\lambda_{i,j}$ are $n$-th roots of unity. We define a map $\varphi : \Z_n \to \Z_n$ by 
\[  \lambda_{1,i} = e(-\varphi(i)/n)    \, .   \]
In particular $\varphi(0) = \varphi(1) = 0 \! \mod n$.
Then
\[  c(i,k) = \varphi(i) + \varphi(k) - \varphi(i+k)  \, , \]
i.e.\ $c$ is a $2$-coboundary arising from $ \varphi$.
We show that $i \varphi(i) = 0 \! \mod n$. 
The above equation implies
\[  \lambda_{i+1,j} 
          = e(-\varphi(j)/n) e( j(\varphi(i) - \varphi(i+1) )/n) \lambda_{i,j}   \]   
so that by induction
\[  \lambda_{i,j} = e(-(i \varphi(j) +  j \varphi(i))/n)  \, . \]
Hence $\lambda_{i,j} =  \lambda_{i,1}^j  \lambda_{j,1}^i$. 
Define $\xi_i = \lambda_{1,i}^i = e(-i\varphi(i)/n)$. Then $\xi_i^2 = \lambda_{i,i} = 1$ so that $\xi_i = \pm 1$. Now
\begin{align*}
\lambda_{i,k} 
&= n e((il+jk)/n) \mathcal{S}_{(i,j),(k,l)} 
= e((il+jk)/n) e(-b_{\Delta}((i,j),(k,l))) \\[1mm] 
&= e(-(i+k)c(i,k)/n) 
= \lambda_{i,i+k} \lambda_{k,i+k}/\lambda_{i+k,i+k} = \lambda_{i,i+k} \lambda_{k,i+k} 
\end{align*}
implies $\xi_i \xi_j = \xi_{i+j}$ so that $\xi_i = 1$ for all $i$ because $\xi_1 = \lambda_{1,1} = 1$. We define representations ${\tilde \phi}_i = e(\varphi(i)/n) \phi_i$ and the corresponding eigenspaces
\[ \widetilde{W}^{(i,j)} 
     = \{ v \in V(g^i) \, | \, {\tilde \phi}_i(g) v = e(j/n) v \}  
     = W^{(i,j - \varphi(i))}  \, . \]
Then
\begin{align*} 
\widetilde{W}^{(i,j)} \boxtimes \widetilde{W}^{(k,l)} 
&= W^{(i,j-\varphi(i)} \boxtimes W^{(k,l-\varphi(k))}   
\cong W^{(i+k,j+l+c(i,k)-\varphi(i)-\varphi(k))} \\
&= W^{(i+k,j+l-\varphi(i+k))} 
= \widetilde{W}^{(i+k,j+l)}
\end{align*}
and
\[  \rho(\widetilde{W}^{(i,j)}) =  \rho(W^{(i,j-\varphi(i))}) = \frac{i(j-\varphi(i))}{n} = \frac{ij}{n} \mod 1  \]
because $i\varphi(i) = 0 \! \mod n$. Furthermore
\[  \widetilde{W}^{(i,j)}{'} = W^{(i,j-\varphi(i))}  \cong W^{(-i,\al(i)-j+\varphi(i))} 
    = W^{(-i,-j-\varphi(-i))} = \widetilde{W}^{(-i,-j)}{'} \]
because $\al(i) = -c(i,-i) = - \varphi(i) - \varphi(-i)\! \mod n$. Finally
\[
\tilde{\mathcal{S}}_{(i,j),(k,l)} 
= \mathcal{S}_{(i,j-\varphi(i)),(k,l-\varphi(k))} 
= \frac{1}{n} e(-b_{\Delta}( (i,j-\varphi(i),(k,l-\varphi(k)) )  
= \frac{1}{n} e(-(il+jk)/n )   \, . 
\]
This finishes the proof. \eop\\[-2mm]

We emphasise that the maps $\phi_i$ chosen in the theorem are not necessarily the same as in Proposition \ref{choicephi}.
  
The fusion group $D$ is given as a set by 
$D = \Z_n \times \Z_n$
with multiplication
$(i,j) + (k,l) = (i+k,j+l+c_d(i,k))$.
This group is isomorphic to the group 
$\Z_{n^2/(n,d)} \times \Z_{(n,d)}$.

\begin{prp}
The discriminant form $D$ is isomorphic to the discriminant form of the even lattice with Gram matrix 
$\left( \begin{smallmatrix} -2t_n & n \\ n & 0 \end{smallmatrix} \right)$.
\end{prp}

For homogeneous $v\in V^G$ the functions $T_{W^{(i,j)}}(v,\tau)$ transform under Zhu's representation $\rho_{V^G}$. Up to a character this representation is the Weil representation $\rho_D$. Since $V$ is holomorphic the central charge $c$ of $V$ is a positive integer multiple of $8$. In particular $c$ is even so that the Weil representation $\rho_D$ is a representation of ${\SL}_2(\Z)$ and 
\[
\rho_{V^G}(M) = \chi_c(M) \rho_D(M)
\]
for all $M \in \text{SL}_2(\Z)$. The character $\chi_c$ is defined by
\[ 
\chi_c(S) = e(c/8) = 1  \quad \text{and} \quad \chi_c(T) = e(-c/24) \, .  \]
The Weil representation is trivial on $\Gamma(N)$ where $N$ is the level of $D$ so that we obtain

\begin{thm} \label{thm.5.14}
Let $V$ be a simple, rational, $C_2$-cofinite, holomorphic vertex operator algebra of CFT-type and central charge $c$. Let $G = \langle g \rangle$ be a finite, cyclic group of automorphisms of $V$ of order $n$. Suppose the fusion group $D$ of $V^G$ has level $N$ under $q_{\Delta}$. Then the trace functions $T_{W^{(i,j)}}(v,\tau)$ and the twisted traces $T(v,i,j,\tau)$ are modular forms of weight $\wt[v]$ for a congruence subgroup of ${\SL}_2(\Z)$ of level $N$ if $24|c$ and level $\text{lcm}(3,N)$ otherwise.
\end{thm}
This result generalises Theorem 1.6 (ii) in \cite{DLM}.

Combining our knowledge of the fusion algebra of $V^G$ with the results about abelian intertwining algebras we can construct new holomorphic vertex operator algebras. 

\begin{thm} \label{orbtheorem}
Let $V$ be a simple, rational, $C_2$-cofinite, holomorphic vertex operator algebra of CFT-type and central charge $c$. Let $G = \langle g \rangle$ be a finite, cyclic group of automorphisms of $V$ of order $n$. Suppose the modules $W^{(i,j)}$, $(i,j) \neq (0,0)$ of $V^G$ have positive conformal weights. Then the direct sum
\[  W = \bigoplus_{i,j \in \Z_n} W^{(i,j)} = \bigoplus_{\gamma\in D} W^{\gamma}   \]
has the structure of an abelian intertwining algebra extending the vertex operator algebra structure of $V^G$ with associated discriminant form $(D,-q_{\Delta})$. Let $H$ be an isotropic subgroup of $D$. Then
\[   W^H = \bigoplus_{\gamma \in H} W^{\gamma}   \]
admits the structure of a simple, rational, $C_2$-cofinite, self-contragredient vertex operator algebra of CFT-type and central charge $c$ extending the vertex operator algebra structure of $V^G$. If $H=H^{\perp}$ then $W^H$ is holomorphic.
\end{thm}

We make some comments on the theorem.

We call $W^H$ the orbifold of $V$ with respect to $G$ and $H$ and denote it by $V^{\orb(G,H)}$.
Any isotropic subgroup $H$ of $D$ of order $n$ satisfies $H = H^{\perp}$.

If $g$ is of type $t \neq  0 \! \mod n$ then $W^H$ can already be obtained from some lower order automorphism with $t = 0 \! \mod n$.

The case when $g$ is of type $t=0 \! \mod n$ is particularly nice. Then it is possible to choose the representations $\phi_i$ such that
\begin{enumerate}[i)]
\item $W^{(i,j)} \boxtimes W^{(k,l)}\cong W^{(i+k,j+l)}$
\item $W^{(i,j)}$ has conformal weight $q_{\Delta}((i,j)) = ij/n \mod 1$
\item ${W^{(i,j)}}' \cong W^{(-i,-j)}$
\item $\mathcal{S}_{(i,j),(k,l)} = \dfrac{1}{n} e(-(jk+il)/n )$, i.e.\ $\sigma(i,k,S) = \lambda_{i,k} = 1$
\item $\mathcal{T}_{(i,j),(k,l)} = e(ij/n-c/24) \delta_{(i,j),(k,l)}$, i.e.\ $\sigma(i,k,T) = e(-c/24)$  
\end{enumerate}
for $i,j,k,l \in \Z_n$. This means that the fusion algebra of $V^G$ is the group algebra of the abelian group $\Z_n \times \Z_n$ with quadratic form
\[  q_{\Delta}((i,j)) = \frac{ij}{n} \mod 1  \, . \]
The direct sum of irreducible $V^G$-modules
\[   W = \bigoplus_{ i,j \in \Z_n} W^{(i,j)}   \]
has the structure of an abelian intertwining algebra and the sum
\[   W^H = \bigoplus_{i \in \Z_n} W^{(i,0)}  \]
over the maximal isotropic subgroup $H = \{ (i,0) \, | \, i \in \Z_n \, \}$ is a simple, rational, $C_2$-cofinite, holomorphic vertex operator algebra of CFT-type which extends the vertex operator algebra structure on $V^G$. In this case we simply write $V^{\orb(G)}$ for the orbifold $W^H$. 

We can define an automorphism $k$ of
$V^{\orb(G)}$
of order $n$ by setting $kv=e(i/n)v$ for $v \in  W^{(i,0)}$. Let $K$ be the cyclic group generated by $k$. 
Then $(V^{\orb(G)})^K = W^{(0,0)} = V^G$ and the twisted modules of $V^{\orb(G)}$ corresponding to $k^j$ are given by $V^{\orb(G)}(k^j) = \bigoplus_{i \in \Z_n} W^{(i,j)}$. Hence
\[ (V^{\orb(G)})^{\orb(K)} = \bigoplus_{j \in \Z_n} W^{(0,j)} = V  \, , \]
i.e.\ there is an inverse orbifold which gives back $V$.

\section{Schellekens' list} \label{Section.schellekens}

In this section we show that Schellekens' classification of $V_1$-structures of meromorphic conformal field theories of central charge $24$ \cite{S} is a rigorous theorem on vertex operator algebras.

\vspace*{2mm}
Let $\mathfrak{g}$ be a finite-dimensional simple Lie algebra over $\C$ and $\mathfrak{h}$ a Cartan subalgebra of $\mathfrak{g}$. We fix a system $\Phi^+$ of positive roots. Let $\theta$ be the corresponding highest root and $\rho$ the Weyl vector, i.e. $\rho = \frac{1}{2} \sum_{\al \in \Phi^+} \al$. We normalise the non-degenerate, invariant, symmetric bilinear form on $\mathfrak{g}$ such that $(\theta,\theta)=2$. This form is related to the Killing form $( \, , \,)_K$ on $\mathfrak{g}$ by $(\ , \,) = \frac{1}{2 h^{\vee}} ( \, , \,)_K$ where $h^{\vee}$ is the dual Coxeter number of $\mathfrak{g}$ (see. \cite{K}, section 6). Note that $\dim(\mathfrak{g}) > (h^{\vee}/2)^2$. 
The affine Lie algebra $\hat{\mathfrak{g}}$ associated to $\mathfrak{g}$ is defined as 
\[ \hat{\mathfrak{g}} = \mathfrak{g} \otimes \C[t,t^{-1}] \oplus \C K  \]
where $K$ is central and 
\[  [a(m),b(n)] = [a,b](m+n) + m \delta_{m+n} (a,b) K  \]
with $a(m) = a \otimes t^m$. Define $ \hat{\mathfrak{g}}_+ = \mathfrak{g} \otimes t\C[t]$, $ \hat{\mathfrak{g}}_- = \mathfrak{g} \otimes t^{-1}\C[t^{-1}]$ and identify $\mathfrak{g}$ with $\mathfrak{g} \otimes 1$. Then 
\[   \hat{\mathfrak{g}} =  \hat{\mathfrak{g}}_+ \oplus  \hat{\mathfrak{g}}_- \oplus \mathfrak{g} \oplus \C K  \, . \]
Let $V$ be a $\mathfrak{g}$-module and $k \in \C$. We can extend the action of $\mathfrak{g}$ to an action of $ \hat{\mathfrak{g}}_+ \oplus \mathfrak{g} \oplus \C K$ by letting $ \hat{\mathfrak{g}}_+$ act trivially and $K$ as $k \, \text{Id}$. Then we form the induced $ \hat{\mathfrak{g}}$-module 
\[  \hat{V}_k = U( \hat{\mathfrak{g}}) \otimes_{U( \hat{\mathfrak{g}}_+ \oplus \mathfrak{g} \oplus \C K)} V  \, . \]
Applied to an irreducible highest-weight module $L(\lambda)$ of $\mathfrak{g}$ we obtain the $ \hat{\mathfrak{g}}$-module $\hat{L}(\lambda)_k$ which we denote by $M_{k,\lambda}$. If $k \neq - h^{\vee}$ then $M_{k,0}$ is a vertex operator algebra of central charge
\[  \frac{k \dim(\mathfrak{g})}{k+h^{\vee}}  \, . \]
The graded pieces of $M_{k,0}$ are $\mathfrak{g}$-modules, and the lowest few are 
\[  M_{k,0} = \C \oplus \mathfrak{g} 
                \oplus ( \mathfrak{g} \oplus \wedge^2(\mathfrak{g}) ) 
                \oplus \ldots   \]
Let $J_{k,\lambda}$ be the maximal proper $\hat{\mathfrak{g}}$-submodule of $M_{k,\lambda}$. Then the quotient $L_{k,\lambda} = M_{k,\lambda}/J_{k,\lambda}$ is an irreducible $ \hat{\mathfrak{g}}$-module. Suppose $k$ is a positive integer. Then $J_{k,0}$ is generated by $e_{\theta}(-1)^{k+1} \vac$ where $e_{\theta}$ is any non-zero element in the root space $\mathfrak{g}_{\theta}$ of $\mathfrak{g}$. By the PBW theorem the lowest graded piece of $J_{k,0}$ has $L_0$-degree $k+1$ and is given by $U(\mathfrak{g}) e_{\theta}(-1)^{k+1} \vac$. Since $(J_{k,0})_{k+1}$ is a finite-dimensional $\mathfrak{g}$-module with highest weight vector $e_{\theta}(-1)^{k+1} \vac$ of weight $(k+1)\theta$ it follows that
\[    (J_{k,0})_{k+1} = L((k+1)\theta)  \]
as a $\mathfrak{g}$-module. Furthermore $L_{k,0}$ is a rational vertex operator algebra whose irreducible modules are the spaces $L_{k,\lambda}$ where $\lambda \in \mathfrak{h}^*$ ranges over the integrable weights satisfying $(\lambda, \theta) \leq k$ (\cite{FZ}, Theorem 3.1.3). The irreducible $L_{k,0}$-module $L_{k,\lambda}$ has conformal weight
\[   \frac{(\lambda, \lambda + 2 \rho)}{2(k+h^{\vee})}   \]
(\cite{K}, Corollary 12.8) and its piece of lowest degree is an irreducible $\mathfrak{g}$-module isomorphic to $L(\lambda)$ because it contains a highest weight vector of weight $\lambda$.

Let $V$ be a simple, rational, $C_2$-cofinite, self-contragredient vertex operator algebra of CFT-type. Then $V$ carries a unique symmetric, invariant, bilinear form $\la \, , \, \ra$ satisfying $\la \vac, \vac \ra = -1$ where $\vac$ denotes the vacuum of $V$ and the subspace $V_1$ of $L_0$-degree $1$ is a Lie algebra under $[a,b] = a_0b$. The bilinear form $\la \, , \, \ra$ is non-degenerate on $V$ and on $V_1$. 

Now we assume in addition that $V$ is holomorphic and of central charge $24$. Then $V_1$ is either trivial or abelian of dimension $24$ or semisimple (see.\ \cite{S} and \cite{DM1}). In the second case $V$ is isomorphic to the vertex algebra of the Leech lattice. We consider the third case where $V_1$ is a semisimple Lie algebra $\mathfrak{g}$. Then the restriction of $\la \, , \, \ra$ to a simple ideal $\mathfrak{g}_i$ of $\mathfrak{g}$ satisfies 
\[  \la \, , \, \ra = k_i ( \, , \, )  \] 
where $( \, , \, )$ is the normalised bilinear form on $\mathfrak{g}_i$ and $k_i$ a positive integer \cite{DM2}. We indicate these integers by writing
\[  \mathfrak{g} = \mathfrak{g}_{1,k_1} \oplus \ldots \oplus \mathfrak{g}_{n,k_n}  \, . \]
We choose a Cartan subalgebra $\mathfrak{h}$ of $\mathfrak{g}$ and decompose $\mathfrak{h}$ accordingly as
\[  \mathfrak{h} = \mathfrak{h}_{1,k_1} \oplus \ldots \oplus \mathfrak{h}_{n,k_n}  \, .  \]  
The map $a(m) \mapsto a_m$ defines a representation of $ \hat{\mathfrak{g}}_{i,k_i}$ on $V$ of level $k_i$. The vertex operator subalgebra of $V$ generated by $V_1$ is isomorphic to 
\[  L_{k_1,0} \otimes \ldots \otimes L_{k_n,0}   \]
and the Virasoro elements of both vertex operator algebras coincide. Since $L_{k_1,0} \otimes \ldots \otimes L_{k_n,0}$ is rational and the $L_0$-eigenspaces of $V$ are finite-dimensional, $V$ decomposes into finitely many irreducible 
$L_{k_1,0} \otimes \ldots \otimes L_{k_n,0}$-modules
\[  V = \bigoplus_{(\lambda_1, \ldots, \lambda_n)} 
         m_{(\lambda_1, \ldots, \lambda_n)}   L_{k_1,\lambda_1} \otimes \ldots \otimes L_{k_n,\lambda_n}  \, . \]
The sum extends over the integrable weights $\lambda_i$ of $\mathfrak{g}_{i,k_i}$ satisfying $(\lambda_i, \theta_i) \leq k_i$ where $\theta_i$ is the highest root of $\mathfrak{g}_{i,k_i}$. The character $\ch_V : H \times \mathfrak{h} \to \C$ defined by
\[  \ch_V (\tau,z) = \tr_V  e^{2 \pi i z} q^{L_0 - 1}   \]
is holomorphic on $H \times \mathfrak{h}$ and transforms as a Jacobi form of weight $0$ and index $1$ (see.\ \cite{G}, section 2, \cite{M0} and \cite{KM}, Theorem 2). In particular it satisfies 
\[  \chi_{V} \left( \frac{a \tau + b}{c\tau + d},\frac{z}{c\tau + d} \right) 
= \exp \left( 2 \pi i \frac{ \la z,z \ra }{2} \frac{ c }{c\tau + d} \right) \chi_{V} (\tau,z)  \]
for all $\left(\begin{smallmatrix} a & b \\ c & d \end{smallmatrix} \right) \in \SL_2(\Z)$. A finite-dimensional $\mathfrak{g}$-module $M$ decomposes into weight spaces and we define for $z \in \mathfrak{h}$ the function
\[  S_M^j(z) = \sum_{\mu \in \Pi(M)} m_{\mu} \mu(z)^j  \]
where $\Pi(M)$ denotes the set of weights of $M$ and $m_{\mu}$ the multiplicity of $\mu$. For example $S_M^0(z) = \dim (M)$.

\begin{thm} \label{scheq}
We have
\[  S_{V_1}^{2}(z) = \frac{1}{12} (\dim(V_1) - 24) \la z,z \ra \\   \]
and for $V_2$
\begin{align*}
S_{V_2}^0(z) =
& \, \, 196884      \\
S_{V_2}^{2}(z) =
& \, \, 32808 \la z,z \ra - 2\dim(V_1) \la z,z \ra \\
S_{V_2}^{4}(z) =
& \, \, 240 \, S_{V_1}^{4}(z) + 15264 \la z,z \ra^2-\dim(V_1) \la z,z \ra^2 \\
S_{V_2}^{6}(z) =
&  \, \, - 504 \, S_{V_1}^{6}(z) + 900 \, S_{V_1}^{4}(z) \la z,z \ra + 11160 \la z,z \ra^3 - 15 \dim(V_1) \la z,z \ra^3 \\
S_{V_2}^{8}(z) =
& \, \, 480 \, S_{V_1}^{8}(z) - 2352 \, S_{V_1}^{6}(z) \la z,z \ra + 2520 \, S_{V_1}^{4}(z) \la z,z \ra^2 + 10920 \la z,z \ra^4 \\
& \, \, - 35 \dim(V_1) \la z,z \ra^4 \\
S_{V_2}^{10}(z) =
& \, \, - 264 \, S_{V_1}^{10}(z) + 2700 \, S_{V_1}^{8}(z) \la z,z \ra - 7560 \, S_{V_1}^{6}(z) \la z,z \ra^2 \\
& \, \, + 6300 \, S_{V_1}^{4}(z)\la z,z \ra^3 + 13230 \la z,z \ra^5 - \frac{315}{4} \dim(V_1) \la z,z \ra^5
\end{align*}
and
\begin{multline*}
48 \, S_{V_2}^{14}(z) - 364 \, S_{V_2}^{12}(z) \la z,z \ra = - 1152 \, S_{V_1}^{14}(z) \\
+ 288288 \, S_{V_1}^{10}(z) \la z,z \ra^2 - 2162160 \, S_{V_1}^{8}(z) \la z,z \ra^3 + 5045040 \, S_{V_1}^{6}(z) \la z,z \ra^4 \\
- 3783780 \, S_{V_1}^{4}(z) \la z,z \ra^5 - 5405400 \la z,z \ra^7 + 45045 \dim(V_1) \la z,z \ra^7.
\end{multline*}
\end{thm}
{\em Proof:}
Let 
\[ P(\tau,z) =  
\exp \left( - (2 \pi i)^2 \frac{ \la z,z \ra }{24} E_2(\tau) \right)
\Delta(\tau)
\ch_V (\tau,z) \]
where 
\[  \Delta(\tau) = q \prod_{m=1}^{\infty} (1-q^m)^{24}  \]
is Dedekind's $\Delta$-function and 
\[  E_2(\tau) = 1 - 24 \sum_{m=1}^{\infty} \sigma(m) q^m  \]
the Eisenstein series of weight $2$. Since $\Delta$ is a modular form of weight $12$ and $E_2$ transforms as
\[  E_2 \left( \frac{a \tau + b}{c\tau + d} \right)  =  
E_2(\tau) (c\tau + d)^2  + \frac{12}{2 \pi i} c (c\tau + d)  \]
for $\left( \begin{smallmatrix} a & b \\ c & d \end{smallmatrix} \right) \in {\SL}_2(\Z)$ we have 
\[  P \left( \frac{a \tau + b}{c\tau + d},\frac{z}{c\tau + d} \right) 
= (c\tau + d)^{12} P(\tau,z)  \, . \]
This implies that the $m$-th coefficient in the Taylor expansion of $P$ in $z$ is a modular form for ${\SL}_2(\Z)$ of weight $12+m$. Using
\begin{align*}
\ch_V (\tau,z) 
& = q^{-1} + \sum_{\lambda \in \Pi(V_1)} m_{\lambda} e^{2 \pi i \lambda(z)} 
    + q \sum_{\lambda \in \Pi(V_2)} m_{\lambda} e^{2 \pi i \lambda(z)} + \ldots \\
& = q^{-1} + \sum_{m=0}^{\infty} \frac{(2 \pi i)^m}{m!} S_{V_1}^{m}(z) 
    + q \sum_{m=0}^{\infty} \frac{(2 \pi i)^m}{m!}S_{V_2}^{m}(z) + \ldots 
\end{align*}
we find that the coefficient of degree $2$ in the Taylor expansion of $P(\tau,z)$ is $(2\pi i)^2$ times
\begin{multline*} 
     - \frac{\la z,z \ra}{24}
     + q   \left( \frac{1}{2} S_{V_1}^2(z) - \frac{\la z,z \ra}{24} (\dim(V_1) - 48) \right)  \\
     + q^2 \left( \frac{1}{2} S_{V_2}^2(z) - 12 S_{V_1}^2(z) - \dim(V_2) \frac{\la z,z \ra}{24} 
                  + 2 \dim(V_1) \la z,z \ra - \frac{63}{2}\la z,z \ra \right)  
     + \, \ldots  \, . 
\end{multline*} 
Since the space of modular forms for ${\SL}_2(\Z)$ of weight $14$ is spanned by 
\[  E_{14}(\tau) = 1 - 24q -196632 q^2 - 38263776 q^3 + \ldots  \]
this implies
\[  \la z,z \ra  = \frac{1}{2} S_{V_1}^2(z) - \frac{\la z,z \ra}{24} (\dim(V_1) - 48) \]
which gives the first equation in the theorem. The second equation is clear. Comparing the coefficients at $q^0$ and $q^2$ we obtain the third equation of the theorem. Looking at higher orders in $z$ and using the relations between the coefficients of modular forms we can derive the remaining relations.  \eop\\[-2mm]

The following result is well-known (see. \cite{S}, \cite{DM1}).
 
\begin{cor} \label{221}
For each simple component $\mathfrak{g}_{i,k_i}$ of $\mathfrak{g}$ we have
\[  \frac{h^{\vee}_i}{k_i} = \frac{\dim(\mathfrak{g}) - 24}{24}   \, . \]
\end{cor}
{\em Proof:}
The first equation of Theorem \ref{scheq} gives
\[  \frac{1}{12} (\dim(\mathfrak{g}) - 24) \la z,z \ra 
                            = \sum_{\mu \in \Phi} \mu(z)^2  \]
where $\Phi$ is the root system of $\mathfrak{g}$. Restricting $z$ to $z_i \in \mathfrak{h}_i$ and using 
\[  \sum_{\mu \in \Phi_i} \mu_i(z_i)^2 =  2 h^{\vee}_i (z_i,z_i)  \]
we obtain the assertion. \eop\\[-2mm]

Now we show how Theorem \ref{scheq} and Corollary \ref{221} can be used to classify the possible affine structures of $V$. 

\begin{prp} 
The equation $h^{\vee}_i/k_i = (\dim(\mathfrak{g}) - 24)/24$ has 221 solutions.
\end{prp}
{\em Proof:}
The equation gives the following inequality
\[  \dim(\mathfrak{g}) 
= \sum_{i=1}^n \dim(\mathfrak{g}_{i,k_i}) 
> \frac{1}{4} \sum_{i=1}^n (h^{\vee}_i)^2 
= \frac{1}{4} \left( \frac{\dim(\mathfrak{g}) - 24}{24} \right)^2 \, \sum_{i=1}^n k_i^2 \, . \]
Since $\sum_{i=1}^n k_i^2   \geq 1$ this implies $\dim(\mathfrak{g}) \leq 2352$. There is a finite set of semisimple Lie algebras satisfying this condition and a computer search yields from this set a list of 221 solutions. \eop \\[-2mm]

Let $V$ be as above with $V_1$ given by one of the solutions
\[  \mathfrak{g} = \mathfrak{g}_{1,k_1} \oplus \ldots \oplus \mathfrak{g}_{n,k_n}  \]
of $h^{\vee}_i/k_i = (\dim(\mathfrak{g}) - 24)/24$. Let $k = (k_1, \ldots, k_n)$, $\lambda = (\lambda_1, \ldots, \lambda_n)$ and $L_{k,\lambda} = L_{k_1,\lambda_1} \otimes \ldots \otimes L_{k_n,\lambda_n}$. Then $V$ decomposes as $L_{k,0}$-module as
\[  V = \bigoplus_{\lambda} m_{\lambda} L_{k,\lambda}  \]
where $\lambda$ ranges over the integrable weights of $\mathfrak{g}$ whose components $\lambda_i$ satisfy $(\lambda_i,\theta_i) \leq k_i$. Since $V_0 = \C \vac$ and $V_1 = \mathfrak{g}$, the modules $L_{k,\lambda}$ with $\lambda \neq 0$ that appear in the decomposition have conformal weight at least $2$. Hence
\[  V_2 = (L_{k,0})_2 \oplus \bigoplus_{\lambda} m_{\lambda} (L_{k,\lambda})_2  \]
where the sum extends over the weights $\lambda$ such that $L_{k,\lambda}$ has conformal weight exactly $2$. The structure of $(L_{k,0})_2$ as $\mathfrak{g}$-module is given by
\[  (L_{k,0})_2 = (M_{k,0})_2 = \mathfrak{g} \oplus \wedge^2(\mathfrak{g})  \]
if $k \geq 2$ and  
\[  (L_{k,0})_2 = (M_{k,0})_2/(J_{k,0})_2 = (\mathfrak{g} \oplus \wedge^2(\mathfrak{g}))/L(2\theta)  \]
with $\theta = (\theta_1, \ldots, \theta_n)$ if $k=1$. If $L_{k,\lambda}$ has conformal weight $2$ then
\[  (L_{k,\lambda})_2 = L(\lambda)  \, . \]
Hence we can decompose the polynomials $S_{V_2}^j$ as
\[  S_{V_2}^j = S_{(L_{k,0})_2}^j + \sum_{\lambda} m_{\lambda} S_{L(\lambda)}^j   \, . \]
Since the polynomials $S_{V_1}^j$ are also known we can write down the second set of equations in Theorem \ref{scheq} explicitly. By equating coefficients we obtain a system of linear equations for the $m_{\lambda}$ which, by the existence of $V$, has a solution in the non-negative integers. We can define an equivalence relation on the weights of $\mathfrak{g}$ by identifying weights which become equal after a permutation of isomorphic components of $\mathfrak{g}$. Then the linear system for the $m_{\lambda}$ gives a reduced system for the multiplicities $m_{[\lambda]} = \sum_{\mu \in [\lambda]} m_{\mu}$ of the classes $[\lambda]$. This system still has a solution in $\Z_{\geq 0}$.

We illustrate this discussion with an example. One of the solutions of the equation $h^{\vee}_i/k_i = (\dim(\mathfrak{g}) - 24)/24$ is $\mathfrak{g} = C_{2,1}^8 A_{2,1}^2$. We write the weights of $C_2 \cong \text{so}_5(\C)$ as $(s_1,s_2)=s_1 \omega_1 + s_2 \omega_2$ where $\omega_1$ denotes the fundamental weight corresponding to the short root of $C_2$ and $\omega_2$ the fundamental weight corresponding to the long root of $C_2$ and similarly for $A_2 \cong \text{sl}_3(\C)$. The irreducible module $L(\mu)$ of $C_2^8 A_2^2$ with highest weight $\mu = ((0, 0)^4, (0, 1)^4, (0, 0)^2)$ has dimension $\dim L(\mu) = 1^4 \cdot 5^4 \cdot 1^2 = 625$. It may be verified that the $L_{k,0}$-module 
$L_{k,\mu}$ has conformal weight $2$. The other irreducible $L_{k,0}$-modules $L_{k,\lambda}$ of conformal weight $2$ can be obtained by permuting the first $8$ components of $\mu$. In particular they all have dimension $625$. Since $\dim (L_{k,0})_2 = 4418$ the second equation of Theorem \ref{scheq} becomes
\[  625 \sum_{\lambda} m_{\lambda} = 192466   \]
where the sum extends over the $\binom{8}{4} = 840$ permutations of the weight $\mu$ introduced above. This equation has no solution in integers so that there is no holomorphic vertex operator algebra of central charge $24$ with $V_1$ isomorphic to $C_{2}^8 A_2^2$.

\begin{thm} \label{Thm.schellekens}
Let $V$ be a simple, rational, $C_2$-cofinite, holomorphic vertex operator algebra of CFT-type and central charge $24$. Then either $V_1 = 0$ or $\dim (V_1) = 24$ and $V$ is isomorphic to the lattice vertex operator algebra of the Leech lattice or $V_1$ is one of $69$ semisimple Lie algebras described in Table~$1$ of \cite{S}.
\end{thm} 
{\em Proof:} If the Lie algebra $V_1$ is non-abelian then it is one of the $221$ solutions $\mathfrak{g}$ of the condition $h^{\vee}_i/k_i = (\dim(\mathfrak{g}) - 24)/24$. For each of these cases we obtain from Theorem \ref{scheq} a system of linear equations on the multiplicities $m_\lambda$ of the decomposition $V_2 = (L_{k,0})_2 \oplus \bigoplus_{\lambda} m_{\lambda} (L_{k,\lambda})_2$. We reduce the system by rewriting it as a system on the multiplicities $m_{[\lambda]} = \sum_{\mu \in [\lambda]} m_\mu$ of equivalence classes as described above. These computations are done using the computer algebra system Magma. For $69$ of the Lie algebras $\mathfrak{g}$ Schellekens \cite{S} gives an explicit candidate decomposition $V = \bigoplus_{\lambda} m_{\lambda} L_k(\lambda)$. We verified in these cases that the multiplicities satisfy the reduced linear system described above. It remains to eliminate the other $\mathfrak{g}$ which we do by ruling out the existence of a solution in $m_{[\lambda]} \in \Z_{\geq 0}$ to the corresponding reduced linear system. Since 
\[   \sum_{\lambda} m_{[\lambda]} \dim(L_k(\lambda)) \leq 196884 \]
and $m_{[\lambda]} \geq 0$ we have $m_{[\lambda]} \leq 196884 / \dim(L_k(\lambda))$. The whole linear system imposes a possibly stronger upper bound on each $m_{[\lambda]}$. We use the simplex algorithm implemented in Mathematica to determine upper bounds for each $m_{[\lambda]}$ subject to the condition $m_{[\mu]} \geq 0$ for all $[\mu] \neq [\lambda]$. For $140$ of the $\mathfrak{g}$ at least one of the upper bounds on the $m_{[\lambda]}$ is strictly negative and the case is thus ruled out. We eliminate a further $4$ cases as we did $C_{2,1}^8A_{2,1}^2$ above by showing that there is no solution in $m_{[\lambda]} \in \Z$. This leaves $8$ cases. From the preceding discussion we have a rational upper bound on each multiplicity $m_{[\lambda]}$. If for some $[\lambda]$ this upper bound is less than $1$ then in fact $m_{[\lambda]}$ should vanish, so we augment the system with the equation $m_{[\lambda]}=0$. Adding this equation to the reduced system we can exclude the cases $A_{3, 16}A_{1, 8}^{5}$, $A_{1, 8}^{2}A_{2, 12}^{3}$, $A_{1, 6}^{6}G_{2, 12}$, $A_{2, 9}^{4}$ and $A_{1, 4}^{2}C_{2, 6}^{3}$. In the cases $A_{1, 16}^{9}$ and $A_{1, 8}^{10}$ we isolate a subset of variables possessing small upper bounds and assign integer values to these variables within their ranges. All possibilities lead to systems that have no solution in $\Q$. The only remaining case is $A_{3,4}^{2}A_{1,2}^{6}$. We supplement the reduced linear system by additional conditions coming from the equality of multiplicities under orbits by fusion with a simple current (see. \cite{S}, Section 3). With these extra equations the system is shown not to possess a solution in $m_{[\lambda]} \in \Z_{\geq 0}$ in the same manner as above. \eop\\[-2mm]

We remark that the first six of the total seven equations in Theorem \ref{scheq} on $V_2$ are not sufficient to rule out $V_1 = \mathfrak{g} = A_{1,8}^{10}$. Only the addition of the seventh equation eliminates this case.

Furthermore entry $62$ in Table 1 of \cite{S} should read $E_{8,2}B_{8,1}$.

\section{Lattice vertex algebras}

In this section we describe some results on lattice vertex algebras and their twisted modules.

Let $L$ be an even lattice and $\varepsilon : L \times L \to \{ \pm 1\}$ a 2-cocycle satisfying $\varepsilon(\al,\al) = (-1)^{(\al,\al)/2}$ and $\varepsilon(\al,\bt)/\varepsilon(\bt,\al)= (-1)^{(\al,\bt)}$. Then $\varepsilon(\al,0) = \varepsilon(0,\al) = 1$, i.e. $\varepsilon$ is normalised. The twisted group algebra $\C[L]_{\varepsilon}$ is the algebra with ba\-sis $\{ e^{\al} | \, \al \in L \}$ and products $e^{\al} e^{\bt}= \varepsilon(\al,\bt) e^{\al + \bt}$. 

Let $g \in \Or(L)$, the orthogonal group of $L$, and $\eta : L \to \{ \pm 1\}$ a function. Then the map $\hat{g}$ on $\C[L]_{\varepsilon}$ defined by $\hat{g}(e^{\al}) = \eta(\al)e^{g(\al)}$ is an automorphism of $\C[L]_{\varepsilon}$ if and only if 
\[  \frac{\eta(\al+\bt)}{\eta(\al)\eta(\bt)} 
      = \frac{\varepsilon(g(\al),g(\bt))}{\varepsilon(\al,\bt)}   \]
for all $\al,\bt \in L$. In this case we call $\hat{g}$ a lift of $g$.

\begin{prp} \label{conjugatelifts}
Let $g \in \Or(L)$ and $g_1, g_2$ be two lifts of $g$ with associated functions $\eta_1, \eta_2 : L \to \{ \pm 1\}$. Suppose $\eta_1 = \eta_2$ on the fixed-point sublattice $L^g$. Then $g_1$ and $g_2$ are conjugate in $\Aut(\C[L]_{\varepsilon})$.
\end{prp}
{\em Proof:} Let $f = 1 - g$ and $M = f(L)$. Set $\eta = \eta_1/\eta_2$. Then $\eta$ is a homomorphism on $L$ which is trivial on $L^g$. Furthermore for all $\al, \bt \in L$ with $f(\al) = f(\bt)$ we have $\eta(\al) = \eta(\bt)$. We construct a homomorphism $\mu : L \to \C^*$ satisfying $\mu \circ f = \eta$ as follows. By the elementary divisor theorem there is a basis $(v_1,\ldots,v_n)$ of $L$ such that $(a_1v_1, \ldots, a_mv_m)$ with $m \leq n$ and non-zero integers $a_i$ is a basis of $M$. For $i = 1, \ldots, m$ choose $\al_i \in L$ such that $a_iv_i = f(\al_i)$ and define $\mu(v_i) \in \C^*$ such that $\mu(v_i)^{a_i} = \eta(\al_i)$. For $i = m+1, \ldots, n$ define $\mu(v_i) \in \C^*$ arbitrarily. Then extend $\mu$ to $L$ by multiplicativity. Note that $\mu(\al)/\mu(g(\al)) = \eta(\al)$ for all $\al \in L$. This implies that the function $h$ on $\C[L]_{\varepsilon}$ defined by $h(e^{\al}) = \mu(\al)e^{\al}$ is an automorphism of $\C[L]_{\varepsilon}$ satisfying $hg_1 = g_2h$.  \eop

\begin{prp} 
Let $g \in \Or(L)$ and suppose $\eta: L \to \{ \pm 1\}$ satisfies the above lifting condition. Then for all $\al \in L^{g^k}$ 
\[ \eta(\al + g(\al) + \ldots + g^{k-1}(\al)) 
   = \eta(\al) \eta(g(\al)) \ldots \eta(g^{k-1}(\al)) \]
if $k$ is odd and 
\[ \eta(\al + g(\al) + \ldots + g^{k-1}(\al)) 
   = (-1)^{(\al,g^{k/2}(\al))} \eta(\al) \eta(g(\al)) \ldots \eta(g^{k-1}(\al)) \] 
if $k$ is even. 
\end{prp}
{\em Proof:} We have 
$\eta(\al+\bt)=\eta(\al)\eta(\bt)\varepsilon(g(\al),g(\bt))/\varepsilon(\al,\bt)$ for all $\al,\bt \in L$. This implies 
\[  \eta(\al+g(\al)) = 
    \eta(\al)\eta(g(\al))\varepsilon(g(\al),g^2(\al))/\varepsilon(\al,g(\al)) \, . \]
More generally, using the $2$-cocycle property of $\varepsilon$, induction on $k$ shows that
\begin{multline*}
\eta(\al + g(\al) + \ldots + g^{k-1}(\al)) = 
\eta(\al) \eta(g(\al)) \ldots \eta(g^{k-1}(\al)) \\
\varepsilon(g(\al) + \ldots + g^{k-1}(\al),g^k(\al))/
\varepsilon(\al,g(\al) + \ldots + g^{k-1}(\al))  \, . 
\end{multline*}
Hence for all $\al \in L^{g^k}$  
\[  \eta(\al + g(\al) + \ldots + g^{k-1}(\al)) = 
        \eta(\al) \eta(g(\al)) \ldots \eta(g^{k-1}(\al)) 
        (-1)^{(\al,g(\al) + \ldots + g^{k-1}(\al)) }  \, .  \]
From this the statement follows easily. \eop\\[-2mm]

There is a lift $\hat{g}$ of $g$ such that $\eta = 1$ on the fixed-point sublattice $L^g$ of $L$. We call such a lift a standard lift of $g$. 

\begin{prp}
Suppose $g$ has odd order $n$ and let $\hat{g}$ be a standard lift of $g$. Then 
\[ \hat{g}^k(e^{\al}) = e^{\al} \]
for all $\al \in L^{g^k}$. In particular $\hat{g}$ has order $n$.
\end{prp}
{\em Proof:} First suppose that $k$ is odd. Then
\[ \hat{g}^k(e^{\al}) 
                  = \eta(\al) \eta(g(\al)) \ldots \eta(g^{k-1}(\al)) e^{\al} 
                  = \eta(\al + g(\al) + \ldots + g^{k-1}(\al)) e^{\al} 
                  = e^{\al}  \]
because $\al + g(\al) + \ldots + g^{k-1}(\al)$ is in $L^g$ and $\eta=1$ on $L^g$. Now suppose that $k$ is even. Then $k+n$ is odd and 
  $ e^{\al} = \hat{g}^{k+n}(e^{\al}) = \hat{g}^{k}(\hat{g}^{n}(e^{\al})) 
           = \hat{g}^{k}(e^{\al}) $
because $\hat{g}$ has order $n$. \eop\\[-2mm]

In the same way one shows 

\begin{prp}
Suppose $g$ has even order $n$ and let $\hat{g}$ be a standard lift of $g$. Then for all $\al \in L^{g^k}$
\[ \hat{g}^k(e^{\al}) = e^{\al} \]
if $k$ is odd and
\[ \hat{g}^k(e^{\al}) = (-1)^{(\al,g^{k/2}(\al))} e^{\al} \]
if $k$ is even. In particular $\hat{g}$ has order $n$ if $(\al,g^{n/2}(\al))$ is even for all $\al \in L$ and order $2n$ otherwise.
\end{prp}

Let $h = L \otimes \C$ and $\hat{h}$ be the Heisenberg algebra corresponding to $h$. Then the group $\Or(L)$ acts on $S(\hat{h}^-)$. The elements of $\Or(L)$ can be lifted as described above to automorphisms of $\C[L]_{\varepsilon}$ and these liftings act naturally on the vertex algebra $V_L = S(\hat{h}^-) \otimes \C[L]_{\varepsilon}$ corresponding to $L$. Proposition \ref{conjugatelifts} implies that two standard lifts of $g \in \Or(L)$ are conjugate in $\Aut(V_L)$. 

Now let $L$ be a positive definite even lattice and $g \in \Or(L)$ of order $n$. Then for even $k$ the map $\xi : L \to \{ \pm 1\}$ defined by $\xi(\al) = (-1)^{(\al,g^{k/2}(\al))}$ is a group homomorphism on $L^{g^k}$ and we define the twisted theta function
\[  \theta_{\xi , L^{g^k} } = \sum_{\al \in L^{g^k}}  \xi(\al) q^{\al^2/2}  \, . \]
Let $\hat{g}$ be a standard lift of $g$. We describe the trace of $\hat{g}^k$ on $V_L$. 

\begin{prp} \label{supertheta}
Suppose $n$ is odd. Then
\[   \tr_{V_L} \hat{g}^k q^{L_0-c/24} = \frac{\theta_{L^{g^k}}(\tau)}{\eta_{g^k}(\tau)}  \]
for all $k$. If $n$ is even then
\[   \tr_{V_L} \hat{g}^k q^{L_0-c/24} = \frac{\theta_{L^{g^k}}(\tau)}{\eta_{g^k}(\tau)}  \]
for odd $k$ and 
\[   \tr_{V_L} \hat{g}^k q^{L_0-c/24} 
                              = \frac{\theta_{\xi, L^{g^k}}(\tau)}{\eta_{g^k}(\tau)}  \]
for even $k$. Here $\eta_{g^k}$ denotes the eta-product corresponding to $g^k$.
\end{prp}

Finally we describe the twisted modules of lattice vertex algebras. We assume for simplicity that $L$ is an even unimodular lattice. Let $g \in \Or(L)$ and $L^{g \perp}$ the orthogonal complement of $L^g$ in $L$. Then the orthogonal projection $\pi : h \to h$ of $h$ onto $L^g \otimes \C$ sends $L$ to the dual lattice ${L^g}'$ of $L^g$. Let $\eta$ be a function on $L$ satisfying the lifting condition. Then $\eta$ is a homomorphism on $L^g$ and there is an element $s_{\eta} \in \pi(h)$ such that 
$\eta(\al) = e((s_{\eta},\al))$ for all $\al \in L^g$. The element $2s_{\eta}$ is in ${L^g}'$ and $s_{\eta}$ is unique up to ${L^g}'$. Let $\hat{g}$ be the lift of $g$ corresponding to $\eta$. A minor variation on the arguments in \cite{BK04} and \cite{DL2} which deal with the case of standard lifts only gives

\begin{thm} \label{confweighttwistedmoduleslatticevoa}
The unique irreducible $\hat{g}$-twisted $V_L$-module is isomorphic as a vector space to
\[  S(\hat{h}_g^-) \otimes e^{s_{\eta}} \C[ \pi(L)] \otimes X  \]
where $X$ is a complex vector space of dimension $d$ with $d^2 = |L^{g \perp}/(1-g)L|$ and $S(\hat{h}_g^-)$ the twisted Fock space. Under this identification the $L_0$-eigenvalue of the vector
\[  \big( a_1(-n_1) \ldots a_l(-n_l) \big) \otimes e^{s_{\eta} + \pi(\al)} \otimes x
\]
is equal to $\sum n_i + (s_{\eta} + \pi(\al))^2/2 + \rho$ where 
\[  \rho = \frac{1}{4n^2} \sum_{j=1}^{n - 1} j(n-j) \dim(h_j) \, . \]
Here $h_j$ denotes the $e(j/n)$-eigenspace of $g$ in $h$.
\end{thm}
 
We need the description of twisted modules for non-standard lifts because sometimes $\hat{g}^m$ is not a standard lift of $g^m$ if $\hat{g}$ is a standard lift of $g$. For example in the case of the order $4$ orbifold in the next section $\hat{g}^2$, where $\hat{g}$ is a standard lift of $g$, is not a standard lift of $g^2$.

\section{Construction of some new holomorphic vertex operator algebras}\label{Constructions}

In this section we construct 5 new holomorphic vertex operator algebras of central charge $24$ as orbifolds of lattice vertex algebras.

We proceed as follows. Let $N(\Phi)$ be a Niemeier lattice with root system $\Phi$ and $g$ an automorphism of $N(\Phi)$ of order $n$. We take a standard lift of $g$ to the vertex algebra $V$ of $N(\Phi)$ which we also denote by $g$. Then the twisted modules $V(g^j)$ have positive conformal weights for $j \neq 0 \! \mod n$ and we can compute the type of $g$ using Theorem \ref{confweighttwistedmoduleslatticevoa}. Suppose that $g$ has type $0$. Then the orbifold of $V$ corresponding to $g$ is given by 
\[   V^{\orb(G)} = \bigoplus_{i \in \Z_n} W^{(i,0)}  \, . \]
The twisted traces $T(\vac,0,j,\tau)$ are described in Proposition \ref{supertheta}. We can determine the other twisted traces $T(\vac,i,j,\tau)$ using the twisted modular invariance  
\[  T(\vac,i,j,M \tau) = T(\vac,(i,j)M ,\tau)  \]
which holds for all $M \in \SL_2(\Z)$. We have 
\[   \ch_{W^{(i,j)}}(\tau) = T_{W^{(i,j)}}(\vac,\tau) 
             = \frac{1}{n} \sum_{l \in \Z_n} e(-lj/n) T(1,i,l,\tau)  \]
so that
\[  \ch_{W^{(i,0)}}(\tau) =  T_{W^{(i,0)}}(\vac,\tau) 
                        = \frac{1}{n} \sum_{j \in \Z_n} T(\vac,i,j,\tau)      \]
and finally
\[  \ch_{V^{\orb(G)}}(\tau) = \sum_{i \in \Z_n}\ch_{W^{(i,0)}}(\tau) = \frac{1}{n} \sum_{i,j \in \Z_n} T(\vac,i,j,\tau)  \, . \]
The constant coefficient in the Fourier expansion of $\ch_{V^{\orb(G)}}$ is $\dim(V^{\orb(G)}_1)$. Once we we have determined the dimension of $V^{\orb(G)}_1$, Theorem \ref{Thm.schellekens} gives the possible Lie algebra structures of $V^{\orb(G)}_1$. We can further restrict the structure of $V^{\orb(G)}_1$ by an argument which is due to Montague \cite{Mon98}. We have seen that there is an automorphism $k$ of order $n$ acting on $V^{\orb(G)}$ such that the corresponding orbifold $(V^{\orb(G)})^{\orb(K)}$ gives back $V$. The fixed-point subalgebra is $(V^{\orb(G)})^K = V^G$ so that $V^G_1$ is the fixed-point subalgebra of an automorphism of $V^{\orb(G)}_1$ whose order divides $n$. Kac classified the finite order automorphisms of finite-dimensional simple Lie algebras and their fixed-point subalgebras (see Theorem 8.6 and Proposition 8.6 in \cite{K}). This reduces the possible affine structures of $V^{\orb(G)}_1$. In two cases we need additional arguments to completely fix this structure.

\subsection*{The affine structure $A_{4,5}^2$ as a $\Z_5$-orbifold}

The lattice $A_4$ can be written as 
$A_4 = \{ \, (x_1,\ldots,x_5) \in \Z^5 \, | \, x_1 + \ldots + x_5 = 0 \, \} \subset \R^5$.
The dual lattice is given by 
$A_4' = \bigcup_{i=0}^4 \, (  (i) +  A_4 )$ 
where $(i) = (\frac{i}{5},\ldots,\frac{i}{5},-\frac{j}{5},\ldots,-\frac{j}{5})$ with $j$ components equal to $\frac{i}{5}$ and $i+j=5$. The lattice $L=A_4^6$ has an automorphism $g$ of order $5$ obtained by acting with a permutation of order $5$ on the coordinates of the first $A_4$-component and a permutation of order $5$ of the remaining $A_4$-components. The characteristic polynomial of $g$ is $(x-1)^{-1}(x^5-1)^5$, i.e.\ $g$ has cycle shape $1^{-1}5^5$. Let $H$ be the isotropic subgroup of $L'/L$ generated by the glue vectors $[1(01441)]$ (see \cite{CS}, Chapter 16). Then the lattice 
$N(A_4^6) = \bigcup_{\gamma \in H} (  \gamma + L )$
is a Niemeier lattice with root system $A_4^6$ and $g$ defines an automorphism of $N(A_4^6)$. The fixed-point sublattice of $g$ in $N(A_4^6)$ is isomorphic to $A_4'(5)$. 
Let $V$ be the vertex operator algebra corresponding to $N(A_4^6)$. We take a standard lift of $g$ which we also denote by $g$. Then $g$ has order $5$ and type $0$. 
The functions 
\begin{align*}
f_1(\tau) 
&= \frac{\theta_{N(A_4^6)}(\tau)}{\Delta(\tau)} = j(\tau) - 600 
= q^{-1} + 144 + 196884q + 21493760q^2 + \ldots \\
f_5(\tau) 
&= \frac{\theta_{N(A_4^6)^g}(\tau)}{\eta_g(\tau)} = 5 + \frac{\eta(\tau)^6}{\eta(5\tau)^6}
= q^{-1} - 1 + 9q + 10q^2 - 30q^3 + 6q^4 + \ldots 
\end{align*}
are modular forms of weight $0$ for $\SL_2(\Z)$ and $\Gamma_0(5)$, respectively and 
\begin{align*}
T(\vac,0,0,\tau) &= f_1(\tau)  \\
T(\vac,0,i,\tau) &= f_5(\tau)  
\end{align*}
for $i \neq 0 \! \mod 5$. An expansion of $f_5$ at the cusp $1/1 = 0$ is given by
\begin{align*}
f_{5,1/1} (\tau) 
&= f_5(S\tau) = 5 + 5^3 \frac{\eta(\tau)^6}{\eta(\tau/5)^6} 
= 5 + 125q^{1/5} + 750q^{2/5} + 3375q^{3/5} + \ldots \\[1mm]
&= f_{5,1/1,0}(\tau) + \ldots + f_{5,1/1,4}(\tau) 
\end{align*}
with $f_{5,1/1,j}(T\tau) = e(j/5)f_{5,1/1,j}(\tau)$. We obtain
\[  \ch_{W^{(0,0)}}(\tau) 
= \frac{1}{5} \sum_{i \in \Z_5} T(\vac,0,i,\tau) = \frac{1}{5} ( f_1 + 4f_5 ) 
= q^{-1} + 28 + 39384q + 4298760q^2 + \ldots  \]  
and 
\begin{align*}
\ch_{W^{(i,0)}}(\tau) 
&= \frac{1}{5} \sum_{j \in \Z_5} T(\vac,i,j,\tau) = \frac{1}{5} \sum_{j \in \Z_5} T(\vac,(i,0)T^j,\tau) 
\end{align*}
\begin{align*}
&= \frac{1}{5} \sum_{j \in \Z_5} T(\vac,(0,i)ST^j,\tau) = \frac{1}{5} \sum_{j \in \Z_5} f_{5,1/1}(T^j\tau) \\[1mm]
&= \frac{1}{5} \sum_{j \in \Z_5} T(\vac,(0,i)ST^j,\tau) = \frac{1}{5} \sum_{j \in \Z_5} f_{5,1/1}(T^j\tau) \\
&= f_{5,1/1,0}(\tau) 
= 5 + 39375q + 4298750q^2 + 172860000q^3 + \ldots
\end{align*}
for $i \neq 0 \! \mod 5$ so that
\[  \ch_{V^{\orb(G)}}(\tau) = \sum_{i \in \Z_5} \ch_{W^{(i,0)}}(\tau) 
= q^{-1} + 48 + 196884q + 21493760q^2 + \ldots \]
Hence the subspace $V^{\orb(G)}_1$ of $V^{\orb(G)}$ has dimension $48$ and by Theorem \ref{Thm.schellekens} is isomorphic as a Lie algebra to $A_1^{16}$, $A_2^6$, $A_1 A_3^3$, $A_4^2$, $A_1 A_5 B_2$, $A_1 D_5$ or $A_6$. It is not difficult to see that $W^{(0,0)}_1 = V^G_1$ is isomorphic to $A_4\C^4$ with $\C^4$ coming from the $4$ orbits of $G$ on the $20$ roots of the first $A_4$-component. Hence $V^{\orb(G)}_1$ admits an automorphism of order dividing $5$ whose fixed-point subalgebra contains $A_4$. The simple components in  the above list possessing such an automorphism are $A_4, A_5, D_5$ and $A_6$. The corresponding fixed-point subalgebras are $A_4, A_4 \C, A_4 \C$ and $A_4\C^2$. Hence $V^{\orb(G)}_1 = A_4^2$ or $A_1 A_5 B_2$. The conformal weight of the spaces $W^{(i,0)}$, $i \neq 0 \! \mod 5$ is $1$. By Lemma 2.2.2 in \cite{SS} this implies that $A_4 \subset V^G_1$ is not only a subalgebra but an ideal in $V^{\orb(G)}_1$. Hence $V^{\orb(G)}$ has the affine structure $A_{4,5}^2$.

\subsection*{The affine structure $C_{4,10}$ as a $\Z_{10}$-orbifold}

Again we consider the Niemeier lattice $N(A_4^6)$ with root system $A_4^6$. Let $g$ be the product of the automorphism of the previous example with $-1$. Then $g$ has order $10$ and cycle shape $1^12^{-1}5^{-5}10^5$. The inner product $(\al,g^k(\al))$ is even for all $\al \in N(A_4^6)^{g^{2k}}$. We take a standard lift of $g$ to the vertex algebra $V$ of $N(A_4^6)$ which we also denote by $g$. Then $g$ has order $10$ and type $0$. 
We define functions
\begin{align*}
f_1(\tau) 
&= j(\tau) - 600 
= q^{-1} + 144 + 196884q + 21493760q^2 + \ldots \\
f_2(\tau) 
&= \frac{\eta(\tau)^{24}}{\eta(2\tau)^{24}} 
= q^{-1} - 24 + 276q - 2048q^2 + 11202q^3 - 49152q^4 + \ldots \\
f_5(\tau) 
&= 5 + \frac{\eta(\tau)^6}{\eta(5\tau)^6} 
= q^{-1} - 1 + 9q + 10q^2 - 30q^3 + 6q^4 - 25q^5 + \ldots \\
f_{10}(\tau) 
&= \frac{\eta(2\tau)\eta(5\tau)^5}{\eta(\tau)\eta(10\tau)^5} 
= q^{-1} + 1 + q + 2q^2 + 2q^3 - 2q^4 - q^5 - 4q^7 + \ldots 
\end{align*}
Then $f_n$ is a modular form for $\Gamma_0(n)$ with trivial character of weight $0$ and
\[   T(\vac,0,i,\tau) = \tr_V g^i \, q^{L_0-1} = f_{10/(i,10)}(\tau) \, . \]
In order to calculate the characters of the modules $W^{(i,0)}$ we need the expansions of these functions at the different cusps. The group $\Gamma_0(2)$ has $2$ cusps and an expansion of $f_2$ at the cusp $1/1$ is given by
\begin{align*}
f_{2,1/1} (\tau) 
&= f_2(S\tau) = 2^{12} \frac{\eta(\tau)^{24}}{\eta(\tau/2)^{24}} 
= 4096q^{1/2} + 98304q + 1228800q^{3/2} + \ldots \\[1mm]
&= f_{2,1/1,0}(\tau) + f_{2,1/1,1}(\tau)
\end{align*}
with $f_{2,1/1,j}(T\tau) = e(j/2)f_{2,1/1,j}(\tau)$.
The group $\Gamma_0(5)$ has $2$ cusps and an expansion of $f_5$ at the cusp $1/1$ is
\begin{align*}
f_{5,1/1} (\tau) 
&= f_5(S\tau) = 5 + 5^3 \frac{\eta(\tau)^6}{\eta(\tau/5)^6} 
= 5 + 125q^{1/5} + 750q^{2/5} + 3375q^{3/5} + \ldots \\[1mm]
&= f_{5,1/1,0}(\tau) + \ldots + f_{5,1/1,4}(\tau) \, .
\end{align*}
Finally the expansions of $f_{10}$ at the cusps $1/5, 1/2$ and $1/1$ of $\Gamma_0(10)$ are given by 
\begin{align*}
f_{10,1/5} (\tau) 
&= f_{10}(M\tau) = - 4 \frac{\eta(\tau/2)\eta(5\tau)^5}{\eta(\tau)\eta(5\tau/2)^5} 
= -4q^{1/2} + 4q + 4q^2 - 4q^{5/2} - 16q^3 + \ldots \\[1mm]
&= f_{10,1/5,0}(\tau) + f_{10,1/5,1}(\tau) 
\end{align*}
where $M = \left( \begin{smallmatrix} 1 & 1 \\ 5 & 6 \end{smallmatrix} \right)$,
\begin{align*}
f_{10,1/2} (\tau) 
&= f_{10}(N\tau) = -  \frac{\eta(2\tau)\eta(\tau/5)^5}{\eta(\tau)\eta(2\tau/5)^5} 
= -1 + 5q^{1/5} - 10q^{2/5} + 15q^{3/5} - 30q^{4/5} + \ldots \\[1mm]
&=f_{10,1/2,0}(\tau) + \ldots + f_{10,1/2,4}(\tau) 
\end{align*}
with $N = \left( \begin{smallmatrix} 1 & -3 \\ 2 & -5 \end{smallmatrix} \right)$ and
\begin{align*}
f_{10,1/1} (\tau) 
&= f_{10}(S\tau) = 4 \frac{\eta(\tau/2)\eta(\tau/5)^5}{\eta(\tau)\eta(\tau/10)^5} 
= 4 + 20q^{1/10} + 60q^{2/10} + 160q^{3/10} + \ldots \\[1mm]
&=f_{10,1/1,0}(\tau) + \ldots + f_{10,1/1,9}(\tau) \, .  
\end{align*}
It follows  
\[  \ch_{W^{(0,0)}}(\tau) = \frac{1}{10} \sum_{i \in \Z_{10}} T(\vac,0,i,\tau) 
= \frac{1}{10} ( f_1 + f_2 + 4f_5 + 4f_{10} ) 
= q^{-1} + 12 + 19720q + \ldots  \]
and
\begin{align*}
\ch_{W^{(5,0)}}(\tau) 
&= \frac{1}{10} \sum_{j \in \Z_{10}} T(\vac,5,j,\tau) 
= \frac{1}{10} \sum_{j \in \Z_5} ( T(\vac,5,j,\tau) + T(\vac,5,5+j,\tau) )  \\
&= \frac{1}{10} ( T(\vac,5,0,\tau) + T(\vac,5,5,\tau) )  \\
& \qquad + \frac{1}{10} \sum_{\substack{j \in \Z_5 \\[0.4mm] (j,5) = 1 }} 
                             ( T(\vac,5,j,\tau) + T(\vac,5,5+j,\tau) )  \\
&=  \frac{1}{10} ( T(\vac,(0,5)S,\tau) + T(\vac,(0,5)ST,\tau) )  \\
& \qquad + \frac{1}{10} \sum_{\substack{j \in \Z_{10} \\[0.4mm] (j,10) = 1 }}
                             ( T(\vac,(0,j)M,\tau) + T(\vac,(0,j)MT,\tau) )  \\
&= \frac{1}{10} (  f_2 (S\tau) + f_2 (ST\tau) ) + \frac{4}{10} (  f_{10} (M\tau) + f_{10} (MT\tau) ) \\[1mm]
&= \frac{1}{10} (  f_{2,1/1} (\tau) + f_{2,1/1} (T\tau) ) + \frac{4}{10} (  f_{10,1/5} (\tau) + f_{10,1/5} (T\tau) ) \\[1mm]
&= \frac{1}{5} f_{2,1/1,0} (\tau) + \frac{4}{5} f_{10,1/5,0}  
= 19664q + 2149584q^2 + 86428864q^3 + \ldots
\end{align*}
In particular $\dim (W^{(5,0)}_1) = 0$. Suppose $(i,10) = 2$. Choose $k \in \Z_{10}$ such that $2k = i \!\mod 10$ and $(k,10)=1$. Then 
\begin{align*}
\ch_{W^{(i,0)}}(\tau) 
&= \frac{1}{10} \sum_{j \in \Z_{10}} T(\vac,i,j,\tau) 
= \frac{1}{10} \sum_{\substack{j \in \Z_{10} \\[0.4mm] (j,2) = 2 }} T(\vac,i,j,\tau) 
   + \frac{1}{10} \sum_{\substack{j \in \Z_{10} \\[0.4mm] (j,2) = 1 }} T(\vac,i,j,\tau) \\
&= \frac{1}{10} \sum_{j \in \Z_5} T(\vac,(0,i)ST^j,\tau) 
   + \frac{1}{10} \sum_{j \in \Z_5} T(\vac,(0,k)NT^j,\tau) \\
&= \frac{1}{10} \sum_{j \in \Z_5} f_5(ST^j \tau) 
   + \frac{1}{10} \sum_{j \in \Z_5} f_{10}(NT^j \tau) \\
&= \frac{1}{10} \sum_{j \in \Z_5} f_{5,1/1}(T^j \tau) 
   + \frac{1}{10} \sum_{j \in \Z_5} f_{10,1/2}(T^j \tau) \\
&= \frac{1}{2} f_{5,1/1,0} + \frac{1}{2} f_{10,1/2,0} 
= 2 + 19715q + 2149170q^2 + 86431120q^3 + \ldots
\end{align*}
Finally for $(i,10) = 1$ we find
\[  \ch_{W^{(i,0)}}(\tau) = f_{10,1/1,0}(\tau) = 4 + 19660q + 2149580q^2 + 86428880q^3 + \ldots  \,  . \]
Hence the character of $V^{\orb(G)}$ is given by
\[  \ch_{V^{\orb(G)}}(\tau) 
= q^{-1} + 36 + 196884q + 21493760q^2 + 864299970q^3 + \ldots  \, . \]
The space $V^{\orb(G)}_1$ has dimension $36$ and therefore is isomorphic as a Lie algebra to $A_1^{12}$, $A_2 D_4$ or $C_4$. The Lie algebra structure of $V_1^G$ is $B_2\C^2$. This excludes the first possibility. Recall that the automorphism $k$ on $V^{\orb(G)}$ acts as $kv=e(i/10)v$ for $v \in W^{(i,0)}$. Hence the fixed-point subalgebra of $k^2$ on $V^{\orb(G)}$ is 
$(V^{\orb(G)})^{K^2} = W^{(0,0)} \oplus W^{(5,0)}$. Since $\dim (W^{(5,0)}_1) = 0$ this implies that $V^G_1 = W^{(0,0)}_1$ is the fixed-point subalgebra of $V^{\orb(G)}_1$ of an automorphism whose order divides $5$ and not just $10$. It follows $V^{\orb(G)}_1=C_4$.

\subsection*{The affine structure $A_{1,1}C_{5,3}G_{2,2}$ as a $\Z_{6}$-orbifold}

Recall that the lattice $E_6$ has 3 conjugacy classes of automorphisms of order $3$ of cycle shape $1^{-3}3^3$, $3^2$ and $1^33^1$. The lattice $L=E_6^4$ has an automorphism of order $3$ acting by a fixed-point free automorphism of order $3$ on the first $E_6$-component and a permutation of order $3$ of the remaining $E_6$-components. Let $g$ be the product of this automorphism with $-1$. Then $g$ has cycle shape $1^32^{-3}3^{-9}6^9$. Let $H$ be the isotropic subgroup of $L'/L$ generated by the glue vectors $[1(012)]$ (see \cite{CS}, Chapter 16). Then the lattice 
$N(E_6^4) = \bigcup_{\gamma \in H} (  \gamma + L )$ is a Niemeier lattice with root system $E_4^6$ and $g$ defines an automorphism of $N(E_6^4)$ satisfying $(-1)^{(\al,g^k(\al))} = 1$ for all $\al \in N(E_6^4)^{g^{2k}}$. 
Let $V$ be the vertex operator algebra corresponding to $N(E_6^4)$. We take a standard lift of $g$ to $V$ which we also denote by $g$. Then $g$ has order $6$ and type $0$. Let $V^{\orb(G)}$ be the corresponding orbifold. We determine the dimension of $V^{\orb(G)}_1$ as in the previous examples. We find $\dim (V^{\orb(G)}_1) = 72$ so that $V^{\orb(G)}_1$ is isomorphic as a Lie algebra to $A_1^{24}$, $A_1^4 A_3^4$, $A_1^3 A_5 D_4$, $A_1^2 C_3 D_5$, $A_1^3 A_7$, $A_1 C_5 G_2$ or $A_1^2 D_6$. The Lie algebra structure of $V_1^G$ is $A_1A_2C_4\C^1$. The only simple components in the list for $V^{\orb(G)}_1$ which admit an automorphism of order dividing $6$ whose fixed-point subalgebra contains $C_4$ are $A_7$ and $D_5$. This implies that $V^{\orb(G)}_1 = A_1^3 A_7$ or $A_1 C_5 G_2$. In the first case $C_4$ is the full fixed-point subalgebra of $A_7$ so that $A_2$ would have to be a subalgebra of $A_1$ which is impossible. Hence $V^{\orb(G)}_1 = A_1 C_5 G_2$.

\subsection*{The affine structure $A_{2,1}B_{2,1}E_{6,4}$ as a $\Z_4$-orbifold}

The lattice $D_6 = \{ (x_1, \ldots ,x_6) \in \Z^6 \, | \, x_1 + \ldots + x_6 = 0 \! \mod 2 \, \} \subset \R^6$ has an automorphism of order $4$ defined by $(x_1,x_2,x_3,x_4,x_5,x_6) \mapsto (x_1,x_2,x_4,x_3,-x_6,x_5)$. Composing this map with the automorphism $(x,y) \mapsto (-y,x)$ on $A_9^2$ we obtain an automorphism $g$ of the Niemeier lattice $N(A_9^2D_6)$ (see \cite{CS}, Chapter 16) of cycle shape $1^22^{-9}4^{10}$. Here $(-1)^{(\al,g(\al))} = -1$ for some elements $\al \in N(A_9^2D_6)^{g^2}$ while $(-1)^{(\al,g^2(\al))} = 1$ for all $\al \in N(A_9^2D_6)$. We take a standard lift of $g$ to the vertex algebra $V$ of $N(A_9^2D_6)$ which we also denote by $g$. Then $g$ has order $4$ and type $0$. The dimension of $V^{\orb(G)}_1$ is $96$ so that $V^{\orb(G)}_1 = A_2^{12}$, $B_2^4 D_4^2$, $A_2^2A_5^2B_2$, $A_2^2A_8$ or $A_2B_2E_6$. The Lie algebra structure of $V_1^G$ is $A_2B_2D_5\C^1$. Out of the list of simple components of $V^{\orb(G)}_1$ only $E_6$ admits an automorphism of order dividing $4$ whose fixed-point subalgebra contains $D_5$. This implies $V^{\orb(G)}_1 = A_2B_2E_6$.

\subsection*{The affine structure $A_{2,6}D_{4,12}$ as a $\Z_6$-orbifold}

The lattice $A_2$ has a fixed-point free automorphism of order $3$ which we denote by $f$. We define an automorphism of $A_2^{12} = \{ (x_1, \ldots ,x_{12}) \, | \, x_i \in A_2 \, \}$ of order $6$ by composing the maps
\begin{align*}
(x_2,x_3,x_7,x_9,x_{10},x_{12}) & \mapsto (x_3,x_{12},x_9,x_{10},x_2,x_7)  \\
(x_6,x_8,x_{11})               & \mapsto (-x_8,-x_{11},-x_6)             \\
(x_1,x_5)                     & \mapsto (-f(x_5),-f(x_1))               \\
x_4                           & \mapsto -f(x_4) 
\end{align*}
This automorphism has cycle shape $1^12^{-2}3^{-3}6^6$ and defines an automorphism $g$ of the Niemeier lattice $N(A_2^{12})$ (see \cite{CS}, Chapter 16) because it preserves the glue group. Then $(-1)^{(\al,g^k(\al))} = 1$ for all $\al \in N(A_2^{12})^{g^{2k}}$. We take a standard lift of $g$ to the vertex algebra $V$ of $N(A_2^{12})$ which we also denote by $g$. Then $g$ has order $6$ and type $0$. The Lie algebra structure of $V_1^G$ is $A_1 A_2 \C^3$. For the dimension of $V^{\orb(G)}_1$ we find $\dim(V^{\orb(G)}_1)=36$ so that $V^{\orb(G)}_1 = A_1^{12}$, $A_2 D_4$ or $C_4$. This implies $V^{\orb(G)}_1 = A_2 D_4$.

\subsection*{Lattice orbifolds}

We summarise our results in the following theorem.

\begin{thm} \label{NewVOATheorem}
There exist holomorphic vertex operator algebras of central charge $c=24$ with the following $5$ affine structures.

\[
\renewcommand{\arraystretch}{1.2}
\begin{array}{c|c|c|c}
\text{Aff.\ structure} & \text{No.\ in \cite{S}} & \text{Niemeier lat.} & \text{Aut.\ order} \\[0.5mm] \hline 
                      &    &           &    \\[-4mm]
A_{2,1} B_{2,1} E_{6,4}       & 28 & A_9^2 D_6 & 4  \\
A_{4,5}^2                   & 9  & A_4^6     & 5 \\
A_{2,6} D_{4,12}             & 3  & A_2^{12}   & 6 \\
A_{1,1} C_{5,3} G_{2,2}       & 21 & E_6^4     & 6 \\
C_{4,10}                    & 4  & A_4^6     & 10 
\end{array} \]

\end{thm}

Together with some recent results by Lam et al.\  (see. \cite{LS3}, \cite{LS4}, \cite{LL}) this implies that all Lie algebras in Schellekens list are realised as $V_1$-algebras of holomorphic vertex operator algebras of central charge $24$.


\begin{thebibliography}{mmmmx}

\setlength{\itemsep}{-0.8ex} 

\bibitem[BK]{BK04} B.\ Bakalov, V.\ G.\ Kac, {\em Twisted modules over lattice vertex algebras}, Lie theory and its applications in physics V, 3--26, World Sci.\ Publ., River Edge, NJ, 2004
\bibitem[B1]{B1} R.\ E.\ Borcherds, {\em Vertex algebras, Kac-Moody algebras and the Monster}, Proc.\ Nat.\ Acad.\ Sci.\ U.S.A.\ {\bf 83} (1986), 3068--3071 
\bibitem[B2]{B2} R.\ E.\ Borcherds, {\em Monstrous moonshine and monstrous Lie superalgebras}, Invent.\ Math.\ {\bf 109} (1992), 405--444 
\bibitem[B3]{B3} R.\ E.\ Borcherds, {\em Automorphic forms with singularities on Grassmannians}, Invent.\ Math.\ {\bf 132} (1998), 491--562 
\bibitem[C]{C} S.\ Carnahan, {\em Generalized moonshine IV: monstrous Lie algebras}, arXiv:1208.6254
\bibitem[CM]{CM} S.\ Carnahan, M.\ Miyamoto, {\em Regularity of fixed-point vertex operator subalgebras}, arXiv:1603.05645 
\bibitem[CS]{CS} J.\ H.\ Conway, N.\ J.\ A.\ Sloane, {\em Sphere packings, lattices and groups}, Grundlehren der Math.\ Wiss.\ {\bf 290}, 3 ed., Springer, New York, 1999
\bibitem[D]{Don94} C.\ Dong, {\em Representations of the moonshine module vertex operator algebra}, Mathematical aspects of conformal and topological field theories and quantum groups, 
Contemp.\ Math., {\bf 175}, Amer.\ Math.\ Soc., Providence, RI, 1994 
\bibitem[DGNO]{DGNO} V.\ Drinfeld, S.\ Gelaki, D.\ Nikshych, V.\ Ostrik, {\em On braided fusion categories I}, Selecta Math.\ (N.S.) {\bf 16} (2010), 1--119
\bibitem[DJX]{DJX} C.\ Dong, X.\ Jiao, F.\ Xu, {\em Quantum dimensions and quantum Galois theory}, Trans.\ Amer.\ Math.\ Soc.\ {\bf 365} (2013), 6441--6469
\bibitem[DL1]{DL1} C.\ Dong, J.\ Lepowsky, {\em Generalized vertex algebras and relative vertex operators}, Progr. Math. {\bf 112}, Birkhäuser, Boston, MA, 1993
\bibitem[DL2]{DL2} C.\ Dong, J.\ Lepowsky, {\em The algebraic structure of relative twisted vertex operators}, J.\ Pure Appl.\ Algebra {\bf 110} (1996), 259--295
\bibitem[DLM]{DLM} C.\ Dong, H.\ Li, G.\ Mason, {\em Modular-invariance of trace functions in orbifold theory and generalized Moonshine}, Comm.\ Math.\ Phys.\ {\bf 214} (2000), 1--56 
\bibitem[DLN]{DLN} C.\ Dong, X.\ Lin, S.-H.\ Ng, {\em Congruence property in conformal field theory}, Algebra Number Theory {\bf 9} (2015), 2121--2166 
\bibitem[DM1]{DM0} C.\ Dong, G.\ Mason, {\em On quantum Galois theory}, Duke Math.\ J.\ {\bf 86} (1997), 305--321
\bibitem[DM2]{DM1} C.\ Dong, G.\ Mason, {\em  Holomorphic vertex operator algebras of small central charge}, Pacific J.\ Math.\ {\bf 213} (2004), 253--266 
\bibitem[DM3]{DM2} C.\ Dong, G.\ Mason, {\em Rational vertex operator algebras and the effective central charge}, Int.\ Math.\ Res.\ Not.\ {\bf 2004}, 2989--3008 
\bibitem[DRX]{DRX} C.\ Dong, L.\ Ren, F.\ Xu, {\em On Orbifold Theory}, 
Adv.\ Math.\ {\bf 321} (2017), 1--30
\bibitem[FLM]{FLM} I.\ Frenkel, J.\ Lepowsky, A.\ Meurman, {\em Vertex operator algebras and the Monster}, Pure Appl.\ Math.\ {\bf 134}, Academic Press, Boston, MA, 1988
\bibitem[FZ]{FZ} I.\ B.\ Frenkel, Y.\ Zhu, {\em Vertex operator algebras associated to representations of affine and Virasoro algebras}, Duke Math.\ J.\ {\bf 66} (1992), 123--168
\bibitem[G]{G} V.\ Gritsenko, {\em 24 faces of the Borcherds modular form $\Phi_{12}$}, arXiv:1203.6503
\bibitem[H1]{H1} Y.-Z.\ Huang, {\em Generalized rationality and a ``Jacobi identity'' for intertwining operator algebras}, Selecta Math.\ (N.S.) {\bf 6} (2000), 225--267 
\bibitem[H2]{H2} Y.-Z.\ Huang, {\em Differential equations and intertwining operators}, Commun.\ Contemp.\ Math.\ {\bf 7} (2005), 375--400 
\bibitem[H3]{H3} Y.-Z.\ Huang, {\em Vertex operator algebras and the Verlinde conjecture}, Commun.\ Contemp.\ Math.\ {\bf 10} (2008), 103--154
\bibitem[H4]{H4} Y.-Z.\ Huang, {\em Rigidity and modularity of vertex tensor categories}, Commun.\ Contemp.\ Math.\ {\bf 10} (2008), 871--911
\bibitem[HS]{HS} G.\ Höhn, N.\ R.\ Scheithauer, {\em A generalized Kac-Moody algebra of rank 14}, J.\ Algebra {\bf 404} (2014), 222--239
\bibitem[K]{K} V.\ G.\ Kac, {\em Infinite-dimensional Lie algebras}, 3 ed., Cambridge University Press, Cambridge, 1990
\bibitem[KM]{KM} M.\ Krauel, G.\ Mason, {\em Vertex operator algebras and weak Jacobi forms}, Internat.\ J.\ Math.\ {\bf 23} (2012) 
\bibitem[LL]{LL} C.-H.\ Lam, X.\ Lin {\em A Holomorphic vertex operator algebra of central charge 24 with weight one Lie algebra $F_{4,6}A_{2,2}$}, arXiv:1612.08123
\bibitem[LS1]{LS3} C.-H.\ Lam, H.\ Shimakura, {\em Orbifold construction of holomorphic vertex operator algebras associated to inner automorphisms}, Comm.\ Math.\ Phys.\ {\bf 342} (2016), 803--841
\bibitem[LS2]{LS4} C.-H.\ Lam, H.\ Shimakura, {\em A holomorphic vertex operator algebra of central charge 24 whose weight one Lie algebra has type $A_ {6,7}$}, Lett.\ Math.\ Phys.\ {\bf 106} (2016), 1575--1585
\bibitem[LY]{LY} C.-H.\ Lam, H.\ Yamauchi, {\em On the structure of framed vertex operator algebras and their pointwise frame stabilizers}, Comm.\ Math.\ Phys.\ {\bf 277} (2008), 237--285
\bibitem[ML]{ML} S.\ MacLane, {\em Cohomology theory of abelian groups}, Proceedings of the International Congress of Mathematicians, Cambridge, Mass., 1950, vol.\ 2, 8--14, Amer.\ Math.\ Soc., Providence, R.\ I., 1952 
\bibitem[M1]{M0} M.\ Miyamoto, {\em A modular invariance on the theta functions defined on vertex operator algebras}, Duke Math.\ J.\ {\bf 101} (2000), 221--236
\bibitem[M2]{M1} M.\ Miyamoto, {\em A $\Z_3$-orbifold theory of lattice vertex operator algebra and $\Z_3$-orbifold constructions}, Symmetries, integrable systems and representations, Springer Proc.\ Math.\ Stat.\ {\bf 40}, 319--344, Springer, Heidelberg, 2013
\bibitem[M3]{M2} M.\ Miyamoto, {\em $C_2$-cofiniteness of cyclic-orbifold models}, Comm.\ Math.\ Phys.\ {\bf 335} (2015), 1279--1286
\bibitem[MT]{MT} M.\ Miyamoto, K.\ Tanabe, {\em Uniform product of $A_{g,n}(V)$ for an orbifold model $V$ and $G$-twisted Zhu algebra}, J.\ Algebra {\bf 274} (2004), 80--96 
\bibitem[M]{Mon98} P.\ S.\ Montague, {\em Conjectured $\Z_2$-orbifold constructions of self-dual conformal field theories at central charge 24 -- the neighborhood graph}, Lett.\ Math.\ Phys.\ {\bf 44} (1998), 105--120 
\bibitem[N]{N} V.\ V.\ Nikulin, {\em Integral symmetric bilinear forms and some of their geometric applications}, Math.\ USSR Izv.\ {\bf 14} (1980), 103--167
\bibitem[SS]{SS} D.\ Sagaki, H.\ Shimakura, {\em Application of a $\Z_3$-orbifold construction to the lattice vertex operator algebras associated to Niemeier lattices}, Trans.\ Amer.\ Math.\ Soc.\ {\bf 368} (2016), 1621--1646
\bibitem[NRS]{NRS} N.\ R.\ Scheithauer, {\em The Weil representation of $\mathit{SL}_2(\Z)$ and some applications}, Int.\ Math.\ Res.\ Notices {\bf 2009} (2009), no. 8, 1488--1545
\bibitem[S]{S} A.\ N.\ Schellekens, {\em Meromorphic $c \! = \! 24$ conformal field theories}, Comm.\ Math.\ Phys.\ {\bf 153} (1993), 159--185
\bibitem[Y]{Y} H.\ Yamauchi, {\em Module categories of simple current extensions of vertex operator algebras}, J.\ Pure Appl.\ Algebra {\bf 189} (2004), 315--328 
\bibitem[Z]{Z} Y.\ Zhu, {\em Modular invariance of characters of vertex operator algebras}, J.\ Amer.\ Math.\ Soc.\ {\bf 9} (1996), 237--302
\end{thebibliography}
\end{document}